\documentclass[a4paper, 11pt]{article}

\usepackage{fullpage}
\usepackage{amsmath, amsfonts, amssymb}
\usepackage{authblk}
\usepackage[sort&compress, numbers]{natbib}
\usepackage{color}
\usepackage{mathrsfs}
\usepackage{graphicx}
\usepackage{multirow}
\usepackage{subfigure}
\usepackage[small,bf]{caption}

\numberwithin{equation}{section}

\renewcommand{\d}{\mathrm{d}}
\newcommand{\R}{\mathbb{R}}

\newcommand{\td}[2]{\frac{\d #1}{\d #2}}

\newcommand{\pd}[2]{\frac{\partial #1}{\partial #2}}

\renewcommand{\vec}[1]{\boldsymbol{#1}}
\newcommand{\tens}[1]{\boldsymbol{\mathsf{#1}}}

\newcommand{\Nw}{N^w} 
\newcommand{\Jw}{Q^w} 
\newcommand{\muw}{\mu^w} 
\newcommand{\Piw}{\Pi^w} 
\newcommand{\Mw}{M^w} 
\newcommand{\Dw}{D^w} 
\newcommand{\nuw}{\nu^w} 

\newcommand{\Nd}{N^d} 
\newcommand{\Jd}{Q^d} 
\newcommand{\Dd}{D^d} 
\newcommand{\nud}{\nu^d} 

\newcommand{\nw}{n^w} 
\newcommand{\phiw}{\phi^w} 
\newcommand{\phid}{\phi^d} 

\newcommand{\Nwhat}{\hat{N}^w} 
\newcommand{\Ndhat}{\hat{N}^d}
\newcommand{\Jwhat}{\hat{Q}^w}
\newcommand{\Jdhat}{\hat{Q}^d}
\newcommand{\muwhat}{\hat{\mu}^w} 

\newcommand{\G}{\mathcal{G}}
\newcommand{\D}{\mathcal{D}}

\newcommand{\unit}[1]{$^{#1}$}

\newcommand{\etal}{\emph{et al.~}}

\newenvironment{subeq}[1]{
\begin{subequations}
#1
\end{subequations}
}

\graphicspath{{figs/}}
\DeclareGraphicsExtensions{.pdf}

\title{Optimal loading of hydrogel-based drug-delivery systems}

\author[1]{Matthew J.~Penn}
\author[2]{Matthew G.~Hennessy\thanks{matthew.hennessy@bristol.ac.uk}}

\affil[1]{Department of Statistics, University of Oxford, Oxford, OX1 3LB, United Kingdom}
\affil[2]{Department of Engineering Mathematics, University of Bristol, Bristol, BS8 1TW, United Kingdom}

\begin{document}

\maketitle

\begin{abstract}
  Drug-loaded hydrogels provide a means to deliver pharmaceutical agents to
specific sites within the body at a controlled rate. 
The aim of this paper is to understand how controlled drug
release can be achieved by tuning the initial distribution of drug
molecules in a hydrogel. A mathematical model is presented for a spherical
drug-loaded hydrogel. The model captures
the nonlinear elasticity of the polymer network and thermodynamics of swelling.
By assuming that the drug molecules are dilute, the equations
for hydrogel swelling and drug transport partially decouple. A fast
optimisation method is developed to accurately compute the optimal initial
drug concentration by minimising the error between the numerical
drug-release profile and a target profile. By taking the target drug efflux
to be piecewise constant, the optimal initial configuration
consists of a central drug-loaded core with isolated drug packets near the
free boundary of the hydrogel. The optimal initial drug concentration is highly
effective 
at mitigating the burst effect, where a large amount of drug is rapidly
released into the environment. The hydrogel stiffness can be used
to further tune the rate of drug release. Although stiffer gels lead to
less swelling and hence reduce the drug diffusivity,
the drug-release kinetics are faster than for soft gels
due to the decreased distance that drug
molecules must travel to reach the free surface. 
\\[1em]
\end{abstract}


\section{Introduction}

A hydrogel is a two-component system consisting of a deformable polymer network
that is saturated with water.
The hydrophilic nature of the polymers creates an energetic
incentive for water molecules to enter the network via diffusion.
In order for the network to accommodate the volume of the water molecules, the
polymers must stretch. Imbibition of water therefore continues until
the energy cost of elastically deforming the polymer network balances the
energy gain of mixing water and polymer. At equilibrium,
the volume of a swollen hydrogel can be tens or even thousands of
times greater than the
volume of the dry polymer network. The ability to precisely
control the degree of swelling via stimuli such as temperature, pH,
and electric fields has led to hydrogels finding use in a diverse range of
applications~\cite{ahmed2015hydrogel}.

Drug-loaded hydrogels have emerged as important  systems for the
controlled and targetted delivery of pharmaceutical
agents~\cite{qiu2001environment, masteikova2003stimuli}.
Controlled delivery means that drug molecules are released at a
prescribed rate; targetted delivery means that the drug molecules are released
at specific locations within the body. The ability to tune the water
content and stiffness of hydrogels leads to excellent
biocompatibility as they are able to mimic a wide range of biological
tissues~\cite{wichterle1960}.
In addition,
the polymer network provides mechanical and chemical shielding that prevents
the degradation of drug molecules before they are released into the body.

Hydrogels provide a pathway to controlled and targetted drug delivery
through their tunable, 
multi-scale architecture and ability to swell when subjected to a
stimulus~\cite{Li2016}. 
Hydrogels possess a macroscopic length scale associated with their overall
size, which can range from microns to millimetres, 
and a nanometric length scale associated with the mesh size
of the polymer network. Both length scales play key roles in the kinetics
of drug release: the gel size controls the distance that drug molecules
must travel to reach the gel surface and be released into the environment,
whereas the mesh size
controls the rate of drug diffusion through the polymer network. 
If the mesh size is much greater than
the hydrodynamic radius of a drug molecule, then drug diffusion is uninhibited
by the presence of the polymer network. However, as the mesh size approaches
the hydrodynamic radius, drug molecules become increasingly immobilised
by the polymer network. Drug molecules that are larger than the mesh size
are effectively entrapped by the network and diffusion is completely
suppressed. Hydrogel swelling
increases the mesh size of the network, thus mobilising drug molecules and
initiating their release into the surroundings. By programming a hydrogel
to swell in response to specific environmental cues, it is possible to deliver
drug payloads to target sites within the body. For example, environmentally
responsive hydrogels have been used to target tumours~\cite{kanamala2016}
and breast cancer
cells~\cite{shirakura2014} by exploiting local increases in pH and
temperature relative to healthy tissue.

Due to the increasingly widespread use of hydrogel-based drug-delivery systems, there is a need for broadly applicable methods that can sensitively control
drug-release profiles~\cite{knipe2014multi}. Although the various length scales
in a hydrogel can be harnessed to alter the drug-release kinetics,
achieving a desired drug-release profile remains a major challenge. The onset
of swelling can drastically change the time scale of drug diffusion
by simultaneously increasing the drug diffusivity and the distance that
drug molecules must travel to reach the free surface.
Moreover, a common problem with drug-delivery systems is the
so-called 
``burst effect'', where a significant proportion of the drug is released in a short initial time frame~\cite{peers2020chitosan}, a phenomenon that can have potentially dangerous effects~\cite{huang2002minimization}. While many advances
in tunable release kinetics have
been made in recent
years~\cite{kikuchi2017hydrogel}, there is still significant scope for
improvement.

The objective of this work is to employ mathematical modelling to explore
the potential of tuning the drug-release profile by varying the
initial drug concentration in drug-loaded hydrogels.
The mathematical model will utilise the theory of nonlinear elasticity
to capture the large deformations of the polymer network that
occur during swelling and the resulting elastic stresses. The generation
of mechanical stress is a particularly important feature
to resolve as it enhances the transport of water molecules through the hydrogel
via stress-assisted diffusion.
An optimisation theory will be developed for computing the initial distribution
of drug molecules that leads to the best approximation of a
target drug-release profile.
The immobility of the drug when the hydrogel is unswollen
means that a non-uniform initial concentration can be experimentally
achieved in a variety of ways~\cite{lee1984effect} and thus has
significant relevance as a control method.

Extensive research on the mathematical modelling of drug-delivery systems has
led to a plethora of literature which has been reviewed
by Siepmann and Siepmann~\cite{siepmann2008}
and Siepmann and Peppas~\cite{siepmann2012}. Caccavo~\cite{Caccavo2019}
compiled a comprehensive overview of models that have been
specifically developed for hydrogel-based drug-delivery systems; these
include simple empirical expressions for data fitting, detailed physical
models based on continuum mechanics, and statistical and neural-network
models. The idea to control the drug-release kinetics via the initial
drug concentration was first proposed by Lee~\cite{lee1986initial}.
Subsequent developments by Lu \etal\cite{lu1998modeling, lu2000photopolymerized}
involved calculating the optimal initial concentration profile, with drug
concentration modelled by the constant-coefficient diffusion equation. 
Georgiadis and Kostoglou~\cite{georgiadis2001optimization} further extended these works to consider the case of a spatially non-uniform diffusion coefficient
as well as allowing this diffusion coefficient to be a free variable.
However, the models used in these optimisation approaches did not account for
the time-dependent swelling of the hydrogel and its subsequent mechanical
response.

The key novelty of this paper therefore arises from combining optimisation
theory with the use of a fully coupled chemo-mechanical model of a hydrogel
based on nonlinear elasticity. Our results reveal that a piecewise-linear
drug-release
profile is best approximated if the initial drug concentration consists
of a central drug-loaded core and a discrete number of drug ``packets''
near the gel surface, the latter of which are
highly localised regions in the gel that are concentrated in drug molecules
and which are separated by wide drug-free zones. Moreover, we find that
the hydrogel stiffness can be used in tandem with the optimal loading to
further tune the drug-release profile and mitigate the burst effect
over a wide range of dosage intervals.

The paper is organised as follows.
In Sec.~\ref{sec:model}, we present a model of a drug-loaded hydrogel.
In Sec.~\ref{sec:benchmarking}, the equilibrium degree of
swelling is computed and its impact on the drug mobility is assessed.
We also explore the drug-release profiles for a uniform loading of drug
molecules.
A theory for the optimal drug loading is developed in
Sec.~\ref{sec:optimisation} and applied to specific scenarios in
Sec.~\ref{sec:case_studies}. The paper concludes in Sec.~\ref{sec:conclusions}.

\section{Mathematical modelling}
\label{sec:model}

We consider the evolution of a spherical, drug-loaded hydrogel
after it is placed in an aqueous environment. The drug molecules
are assumed to be too large to move through the polymer network when the
gel is in its initial, undeformed state.
For simplicity, it is assumed that the system
remains axisymmetric during swelling and drug release.

Due to the large deformations that occur during swelling, the mechanical
response of the hydrogel is described using the framework of nonlinear
elasticity. On sufficiently long time scales, the polymers may
rearrange to relax the elastic stress,
resulting in a viscoelastic material response~\cite{Caccavo2018, ferreira2018}.
Moreover, the polymer network may degrade~\cite{ashley2013, de2021}.
Neither of these features will be considered here.

Several thermodynamically consistent hydrogel models based
on finite-strain elasticity have been
proposed~\cite{Hong2008, Chester2010, Drozdov2016}.
These models typically assume that the polymer network is swollen by a
solute-free
liquid and describe the evolution of the mixture towards a state of minimum
energy. The addition of a solute, such as
drug molecules,
can alter the energetic landscape of the system
and impact the transport of fluid via cross-diffusion, both of which
require extended models to capture~\cite{Celora2022}.
In the context of drug-delivery systems,
the volume (or mass) fraction of
drug molecules is often small~\cite{Caccavo2018}. Consequently,
the chemo-mechanics of swelling will not be
strongly influenced by the presence of drug molecules~\cite{Volpert2018}
and cross-diffusion can be neglected. In formulating the model below,
we invoke the assumption that the drug molecules are dilute.
As a result, there will be a partial decoupling of the model: governing
equations for the hydrogel can be formulated and solved independently
from those for the drug. This decoupling will be the key to developing
a fast algorithm for optimising the initial drug distribution.

\subsection{Bulk equations for the hydrogel}

The governing equations for the hydrogel have been derived using
thermodynamic arguments by Hennessy \etal\cite{Hennessy2020}; here we
specialise the results to a spherical geometry. The equations 
are formulated in terms of
Lagrangian coordinates $\vec{X}$ associated with the
stress-free reference configuration, which is taken to be a dry
hydrogel with radius $R_0$. The use of Lagrangian coordinates avoids the
introduction of a free boundary into the problem. 
The Lagrangian gradient operator
is denoted by $\nabla_{\vec{X}}$ and is expressed in terms of the usual
spherical coordinates. The Eulerian coordinates
associated with the current (swollen) configuration
are denoted by $\vec{x}$. For the axisymmetric configurations considered
here, we can write $\vec{X} = R \vec{e}_r$ and $\vec{x} = r(R,t) \vec{e}_r$,
where $R$ and $r$ are Lagrangian and Eulerian radial coordinates, respectively,
and $\vec{e}_r$ is the radial basis vector. The deformation gradient tensor
$\tens{F} = \nabla_{\vec{X}} \vec{x}$ describes
the local distortion of material elements,
whereas its determinant,
$J = \det \tens{F}$, describes the volumetric changes of material elements.
For axisymmetric deformations in spherical geometries,
the appropriate form of the deformation gradient tensor is readily calculated
as
\begin{align}
  \tens{F} = \lambda_r \vec{e}_r \otimes \vec{e}_r +
  \lambda_\theta \vec{e}_\theta \otimes \vec{e}_\theta +
  \lambda_\varphi \vec{e}_\varphi \otimes \vec{e}_\varphi,
\end{align}
where $\vec{e}_\theta$ and $\vec{e}_\varphi$ are the polar and azimuthal
basis vectors, $\otimes$ denotes the dyadic product of two vectors, and
\begin{align}
  \lambda_r = \pd{r}{R}, \quad \lambda_{\theta} = \lambda_\varphi = \frac{r}{R},
  \label{eqn:lambda}
\end{align}
are the principal stretches in the radial, polar, and azimuthal directions.
The polymer network and the imbibing fluid, assumed to be water, are treated
as incompressible. This assumption, in combination with the limit of a
dilute drug, implies that volumetric changes in material elements must be
solely associated with the imbibition of water molecules. We therefore
impose a molecular incompressibility condition given by
\begin{align}
  J = \lambda_r \lambda_\theta \lambda_\varphi = 1 + \nuw \Nw,
\end{align}
where $\nuw$ and $\Nw$ are the molecular volume and
nominal concentration of water, respectively.
Nominal concentrations are expressed in terms of the number of molecules
per unit reference (undeformed) volume. The actual concentration
of water, \emph{i.e.}\ the number of molecules per current (deformed) volume,
is defined as $\nw = \Nw / J$. The volume fraction of water $\phiw$ can
then be defined as $\phiw = \nuw \nw$.


The conservation of water can be expressed as
\begin{align}
  \pd{\Nw}{t} + \frac{1}{R^2}\pd{}{R}\left(R^2 \Jw\right) = 0,
\end{align}
where $t$ is time and $\Jw$ is the nominal diffusive flux given by
\begin{align}
  \Jw = -\Mw \pd{\muw}{R}.
\end{align}
The water mobility $\Mw$ is defined as
\begin{align}
  \Mw = \frac{\Nw}{k_B T}\frac{\Dw(J)}{\lambda_r^2},
  \label{eqn:M}
\end{align}
with $k_B$ denoting Boltzmann's constant, $T$ the absolute
temperature, $\Dw$ the diffusivity of water in the polymer network,
and $\muw$ the chemical potential of water.
The factor of $\lambda_r^{-2}$ in \eqref{eqn:M} is a result of mapping
Fick's law in the current configuration to the reference configuration.
The water diffusivity is expressed as
$\Dw = \Dw_0 J^{a}$, where $a$ is a positive parameter that
characterises how strongly the rate of diffusion increases as the
polymer network expands. Typically, $a = 1.5$; see
Bertrand \etal\cite{Bertrand2016}. 

The chemical potential of water can be expressed as
\begin{align}
  \muw = \muw_0 + \nuw (\Piw + p),
\end{align}
where $\muw_0$ is the chemical potential of a pure bath of water,
$\Piw$ is the osmotic pressure of water, and $p$
is the mechanical pressure. The osmotic pressure captures fluid transport
that is driven by concentration gradients and is given by
\begin{align}
  \Piw = \frac{k_B T}{\nuw}\left[\log\left(\frac{\nuw \Nw}{1+\nuw\Nw}\right) + \frac{1 + \nuw \Nw + \chi}{(1 + \nuw \Nw)^2}\right].
\end{align}
The Flory parameter $\chi$ describes the strength of energetically unfavourable
interactions between polymers and water molecules.
Typically, a large value of $\chi$ corresponds to a low degree of swelling,
as it becomes energetically costly for fluid and polymers to mix. The
dependence of the chemical potential on the mechanical pressure captures
transport of fluid down stress gradients and leads to
stress-assisted diffusion.

Conservation of linear momentum is given by
\subeq{
  \label{eqn:stress}
\begin{align}
  \frac{1}{R^2}\pd{}{R}\left(R^2 S_{r}\right) - \frac{1}{R}\left(S_{\theta} + S_{\varphi}\right) = 0, \\
  S_{\theta} - S_{\varphi} = 0,
\end{align}}
where $S_i$ are the principal first Piola--Kirchhoff stresses. The
components of the stress can be decomposed into elastic components
$\Sigma_i$ and a pressure component such that
\begin{align}
  S_i = \Sigma_i - p J \lambda_i^{-1}.
  \label{eqn:S_i}
\end{align} 
The hydrogel is assumed to be a hyperelastic material described
by a neo-Hookean strain energy. Consequently,
the elastic components of the stress
can be written as
\begin{align}
  \Sigma_r = G(\lambda_r - \lambda_r^{-1}), \quad \Sigma_\theta = \Sigma_\varphi = G(\lambda_\theta - \lambda_\theta^{-1}),
\end{align}
where $G$ is the shear modulus of the polymer network.

\subsection{Bulk equations for drug diffusion}

As with the hydrogel model, the equations that govern the transport of
drug molecules are written in terms of Lagrangian coordinates. We let
$\Nd$ represent the nominal concentration of drug molecules, which
must obey the conservation law
\begin{align}
  \pd{\Nd}{t} + \frac{1}{R^2}\pd{}{R}\left(R^2 \Jd\right) = 0.
\end{align}
The volume fraction of drug $\phid$ is defined as
$\phid = \nud \Nd / J$, where $\nud$ is the volume of a
drug molecule. The dilute-drug limit requires $\phid \ll 1$. 
The diffusive flux of drug molecules, $\Jd$, is given by
\begin{align}
  \Jd = -\frac{\Dd(J)}{\lambda_r^2} \pd{\Nd}{R}.
\end{align}
Various forms of the drug diffusivity $\Dd$ appear in the literature.
A common choice is a Fujita-type expression, in which the drug diffusivity
is assumed to exponentially increase with the water concentration. Such
forms are suitable for models that neglect the mechanics
of the polymer network~\cite{caccavo2015} or which only consider small deformations~\cite{ferreira2018} because the water concentration
can serve as a proxy for the degree of swelling that occurs in each material
element. Given that our model explicitly captures finite deformations
of the polymer network, 
we choose an expression for the drug diffusivity
based on free-volume theory~\cite{Caccavo2018}:
\begin{align}
  \Dd(J) = \Dd_\infty \exp\left[-\frac{\beta}{J - 1}\right].
\end{align}
The fitting parameter $\beta > 0$ controls how strongly the drug diffusivity
increases with the volumetric expansion of the polymer network. The
value of $\Dd_\infty$ describes the diffusivity of drug
molecules when they are uninhibited by the polymer network.

\subsection{Boundary and initial conditions}

At the centre of the hydrogel, $R = 0$, we impose
\begin{align}
  r(0,t) = 0,
\end{align}
which ensures that the origin in the current state is mapped to the origin
in the reference state. In addition, we impose no-flux conditions on the
water and drug molecules:
\subeq{
\begin{align}
  \Jw(0,t) &= 0, \\
  \Jd(0,t) &= 0.
\end{align}}
At the free surface of the hydrogel, $R = R_0$,
we impose continuity of the chemical
potential of water and a stress-free condition, which leads to
\subeq{\label{bc:free_surf}
\begin{align}
  \muw(R_0,t) &= \muw_0, \\
  S_r(R_0,t) &= 0. \label{bc:stress}
\end{align}}
In \eqref{bc:free_surf}, we have set the pressure in the surrounding
water to be zero. The surrounding environment is assumed to be
a perfect sink for the drug. Therefore, we impose that the concentration
of drug at the free surface of the hydrogel is zero:
\begin{align}
  \Nd(R_0,t) = 0.
\end{align}
Consequently, the steady-state configuration will correspond to a
swollen, drug-free hydrogel. 


Initially, the hydrogel is in a dry state that does not contain water
molecules but which is loaded with drug molecules. The initial
concentration of drug is denoted by $d$. Therefore, the initial conditions
for the model are
\begin{align}
  \Nw(R,0) = 0, \quad \Nd(R,0) = d(R).
\end{align}
Since the hydrogel is initially undeformed, the nominal and actual
concentrations coincide at $t = 0$.
Thus, the initial volume fraction
of drug is given by $\phid(R,0) = \nud d(R)$. The initial ratio
of the total drug volume to the total volume can be defined as
\begin{align}
  \epsilon = \frac{3}{R_0^3} \int_{0}^{R_0} \nud d(R) R^2\,\d R < 1.
  \label{eqn:d_cons}
\end{align}
We will impose a value of $\epsilon \ll 1$ to
assist in defining the target drug-release profiles.
Given that both $\phid$ and $\epsilon$ represent volume fractions of
drug, we will refer to $\phid$ as the local drug fraction and
$\epsilon$ as the global drug fraction.


\subsection{Drug efflux and target profiles}
\label{sec:flux_target}

The flux of drug molecules out of the hydrogel, hereafter referred to as the
efflux, is defined as
\begin{align}
  F(t) = -\td{}{t}\left[ 4 \pi \int_{0}^{R_0} \Nd(R,t)R^2\, \d R\right]
  = 4 \pi R_0^2 \Jd(R_0,t).
  \label{eqn:F}
\end{align}
By integrating the first equality in \eqref{eqn:F} in time and using
\eqref{eqn:d_cons}, we see that
the efflux $F$ must satisfy the condition
\begin{align}
  \int_{0}^{\infty} F(t)\,\d t = \frac{4\pi R_0^3 \epsilon}{3\nud}.
  \label{eqn:F_cons}
\end{align}
From this point forward, \eqref{eqn:F_cons} will be used in place of
\eqref{eqn:d_cons}.
We now let $A(t)$ denote a target flux profile that
is desirable to achieve in practical situations. We will be
particularly concerned with piecewise-constant target profiles of the form
\subeq{  \label{eqn:A}
\begin{align}
  A(t) = \begin{cases}
    A_0, \quad &0 \leq t \leq \tau, \\
    0, \quad &\text{otherwise},
  \end{cases}
\end{align}
where $\tau$ is referred to as the drug-release period and it describes the
amount of time needed for all of the drug molecules to be released from the
hydrogel. The constant $A_0$ is determined by imposing
the constraint in \eqref{eqn:F_cons}, which ensures that the actual efflux $F$
and the target efflux $A$ lead to the same amount of drug being delivered.
Therefore, we must have that
\begin{align}
  A_0 = \frac{4\pi R_0^3 \epsilon}{3 \tau \nud}.
\end{align}}
The target profile given by \eqref{eqn:A}
has significant physical meaning, as often the goal of drug
delivery through hydrogels is to steadily release the drug over
a set period of time \cite{fu2011experimental}.
The aim of this paper is to determine the initial concentration
of drug molecules in the gel that minimises the error between
the efflux $F(t)$ and the target profile $A(t)$. 
The piecewise-constant target profiles in \eqref{eqn:A} provide an excellent
means of testing the robustness
of our optimisation approach because 
the discontinunity when $t = \tau$ is difficult to approximate.

\subsection{Non-dimensionalisation}

The governing equations are written in dimensionless form using
the initial gel radius $R_0$ as the length scale and $R_0^2 / \Dw_0$ as
the time scale. Thus, we write $R = R_0 \hat{R}$, $r = R_0 \hat{r}$, and
$t = (R_0^2 / \Dw_0) \hat{t}$, where hats are used to denote
non-dimensional quantities.
The chemical potential of water is written as $\muw = \muw_0 + k_B T \muwhat$.
The nominal concentrations, the diffusive fluxes, and the drug efflux
are written as
\begin{align}
  \Nw = \frac{1}{\nuw}\,\Nwhat,
  \quad
  \Nd = \frac{\epsilon}{\nud}\,\Ndhat,
  \quad
  \Jw =  \frac{\Dw_0}{\nuw R_0} \Jwhat,
  \quad 
  \Jd = \frac{\epsilon \Dd_\infty}{\nud R_0} \Jdhat,
  \quad
  F = \frac{\epsilon R_0 \Dw_0}{\nud} \hat{F}.
\end{align}
The elastic stresses and the pressure are non-dimensionalised according to
$\Sigma_r = G \hat{E}_r$, $\Sigma_\theta = G \hat{E}_\theta$, and $p = G \hat{p}$.
In the dimensionless equations presented below, the hats on the variables
will be dropped.

The dimensionless equations for the hydrogel are as follows. The
conservation of water reads as
\subeq{\label{nd:gel:fluid}
\begin{align}
  \pd{\Nw}{t} &= \frac{1}{R^2}\pd{}{R}\left(R^2 \Mw
                \pd{\muw}{R}\right),
                \label{nd:gel:Nw}
\end{align}
where the water mobility is given by $\Mw = \Nw J^{a}/\lambda_r^2$.
The chemical potential of water in the gel is given by
\begin{align}
  \muw = \log\left(\frac{\Nw}{1+\Nw}\right) + \frac{1 + \Nw + \chi}{(1 + \Nw)^2} + \G p,
  \label{nd:gel:muw}
\end{align}}
where $\G = \nuw G / (k_B T)$ is a non-dimensional elastic modulus
that characterises the energy increase due to elastic deformations relative
to the energy decrease of inserting a water molecule into the polymer
network. The radial stress balance can be reduced to
\subeq{\label{nd:gel:mech}
\begin{align}
  \pd{\Sigma_r}{R} + \frac{2(\Sigma_r - \Sigma_\theta)}{R} &= \lambda_{\theta}^2 \pd{p}{R}.
\end{align}
The elastic components of the stress are written as
\begin{align}
  \Sigma_r = \lambda_r - \lambda_r^{-1}, \quad
  \Sigma_\theta = \lambda_\theta - \lambda_\theta^{-1},
\end{align}}
where the non-dimensional expressions for the radial and orthoradial
stretches $\lambda_r$ and $\lambda_\theta$
are identical to those in \eqref{eqn:lambda}.
The molecular incompressibility condition can be formulated as
\begin{align}
  J = \lambda_r\lambda_\theta^2 = 1 + \Nw.
  \label{nd:gel:ic}
\end{align}
Equations \eqref{nd:gel:fluid}--\eqref{nd:gel:ic}
are solved with the following boundary and
initial conditions:
\begin{align}
  r(0,t) = 0,
  \quad
  \left.\pd{\muw}{R}\right|_{R=0} = 0,
  \quad
  \muw(1,t) = 0,
  \quad
  S_r(1,t) = 0,
  \quad
  \Nw(R,0) = 0,
  \label{nd:gel:bc}
\end{align}
where the non-dimensional total radial stress $S_r$ has the same form
as in \eqref{eqn:S_i}.

The diffusion equation that governs the transport of drug molecules
through the hydrogel can be written in dimensionless form as
\subeq{\label{nd:drug:all}
\begin{align}
  \pd{\Nd}{t} &= \frac{\D}{R^2}\pd{}{R}\left(R^2 \frac{\Dd(J)}{\lambda_r^2}
  \pd{\Nd}{R}\right),
  \label{nd:drug:eqn}
\end{align}
where $\D = \Dd_\infty/\Dw_0$. 
The diffusivity of drug molecules is given by
\begin{align}
  \quad
  \Dd = \exp\left(-\frac{\beta}{J - 1}\right).
  \label{eqn:Dd}
\end{align}
The boundary and initial conditions for the
drug concentration are
\begin{align}
  \left.\pd{\Nd}{R}\right|_{r=0} = 0, \quad
  \Nd(1,t) = 0, \quad
  \Nd(R,0) = d(R).
  \label{nd:drug_bc}
\end{align}}
Due to the choice of non-dimensionalisation, the initial volume fraction
of drug is given by $\phid(R,0) = \epsilon d(R)$.
The dilute
limit therefore requires that $d = O(\epsilon^{-1})$ as $\epsilon \to 0$.
The non-dimensional drug efflux is defined as and must satisfy
\begin{align}
  F(t) = -\td{}{t}\left(4 \pi \int_{0}^{1} \Nd(R,t)R^2\, \d R\right),
  \qquad
  \int_{0}^{\infty} F(t)\,\d t = \frac{4}{3}\pi.
\end{align}
Similarly, non-dimensionalising the piecewise-constant target profiles in
\eqref{eqn:A} leads to
\begin{align}
  A(t) = \begin{cases}
    \displaystyle 4 \pi / (3 \tau), \quad &0 \leq t \leq \tau, \\
    0, \quad &\text{otherwise}.
  \end{cases}
               \label{nd:A}
\end{align}

\subsection{Parameter estimation}

The initial radius of the hydrogel is assumed to be 2~mm.
Gels of this size would likely be surgically implanted into the body or
placed directly onto the skin for transdermal drug delivery~\cite{Li2016}.
Water has a molar volume of
$18\cdot 10^{-6}$~m\unit{3}~mol\unit{-1}. Dividing by Avogadro's number gives a molecular volume of $\nuw = 3.0\cdot 10^{-29}$~m\unit{3}.
The shear modulus of the gel, $G$,
is taken to be a control parameter. Typical values range from about 10~kPa
to 1000~kPa. The Flory interaction parameter $\chi$ depends on the specific
type of polymers used to create the gel and is generally
a function of composition.
However, its value often lies between 0 (for athermal mixtures) and 3. For
simplicity, we treat $\chi$ as a constant. Drozdov \etal\cite{Drozdov2016}
fitted a similar hydrogel model to experimental data and reported a solvent
diffusion coefficient of $\Dw_0 \sim 10^{-9}$~m\unit{2}~s\unit{-1}. We follow
Caccavo \etal\cite{Caccavo2018} and assume that $\beta = 1.0$ and
that the drug diffusivity  $\Dd_\infty$ lies in the range
$10^{-12}$~m\unit{2}~s\unit{-1} to
$10^{-10}$~m\unit{2}~s\unit{-1}. For reference, the diffusivity of
paracetamol in water is roughly
$6.5 \cdot 10^{-10}$ m\unit{2}~s\unit{-1} \cite{ribeiro2012}. The diffusivity
of larger macromolecules such as proteins is expected to be substantially
smaller. 
The temperature is fixed at 293~K.

Using these parameter estimates, we find that the non-dimensional shear
modulus, $\G$, is between $7\cdot 10^{-5}$ and $7 \cdot  10^{-3}$,
the smallness of which is characteristic of a soft solid. The diffusivity
ratio $\D$ lies between $10^{-3}$ to $10^{-1}$. The time scale of
fluid diffusion, $R_0^2 / \Dw_0$, is roughly one hour. Therefore, we
assume that the non-dimensional drug-release period $\tau$ appearing in the
target flux profile \eqref{eqn:A} ranges from 6 to 24, corresponding to
drug release over a 6- to 24-hour window.


\subsection{Finite-difference discretisation}
\label{sec:numerics}

The non-dimensional hydrogel model \eqref{nd:gel:fluid}--\eqref{nd:gel:bc}
is solved using a semi-implicit
finite-difference method with a staggered grid.
The Lagrangian spatial domain, $0 \leq R \leq 1$,
is discretised into cells of uniform width. The Eulerian radial coordinate $r$ is solved for
on cell edges whereas the pressure $p$ and the
nominal fluid fraction $\Nw$ are solved for on cell midpoints.
The conservation equation for the fluid \eqref{nd:gel:Nw} is discretised in time
using Euler's method. All quantities in the dimensionless system
\eqref{nd:gel:fluid}--\eqref{nd:gel:bc}
are treated implicitly, with the exception of the mobility $\Mw$
in \eqref{nd:gel:Nw},
which is treated explicitly. This choice provides greater numerical stability
during the first few time steps, where large gradients in the fluid fraction
and radial stretch develop. The resulting nonlinear algebraic system is solved
using Newton's method at each time step.

Once the solution to the gel problem is obtained, the linear diffusion problem for the drug \eqref{nd:drug:all} is solved using an implicit Euler method.
The equations are discretised using the same staggered grid for the hydrogel, with the drug concentration $\Nd$ found on cell midpoints.

\section{Benchmarking}
\label{sec:benchmarking}

The equilibrium states provide valuable information about how swollen the
gel becomes for a given set of parameters. From this information, it is possible
to assess the change in drug diffusivity that occurs during the swelling
process. 
The equilibrium states correspond to homogeneous gels with swelling
ratio $J = J_\infty$.
The radial and orthoradial stretches are therefore equal and given by
$\lambda_r = \lambda_\theta = J_\infty^{1/3}$. The chemical potential of water
$\muw$ and the radial component of the first Piola--Kirchhoff stress $S_r$
are both uniform and, from the boundary conditions in \eqref{nd:gel:bc},
equal to zero. The latter can be used to obtain an expression for
the pressure given by $p = J_\infty^{-1/3} - J_\infty^{-1}$. 
Thus, by writing
the nominal drug concentration in terms of the swelling ratio $J$
using the incompressibility condition \eqref{nd:gel:ic}
and eliminating the pressure in the
chemical potential \eqref{nd:gel:muw},
we find that the equilibrium states satisfy 
\begin{align}
  \log(1 - J_\infty^{-1}) + J_\infty^{-1} + \chi J_\infty^{-2} + \G(J_\infty^{-1/3} - J_\infty^{-1}) = 0.
  \label{eqn:eq}
\end{align}
Once the unique solution for $J_\infty$ is obtained from \eqref{eqn:eq}, the
equilibrium drug diffusivity $\Dd_\infty = \Dd(J_\infty)$, which
represents the maximum value of the drug diffusivity during swelling,
can be computed by evaluating \eqref{eqn:Dd}.

\begin{figure}
  \centering
  \includegraphics[width=0.32\textwidth]{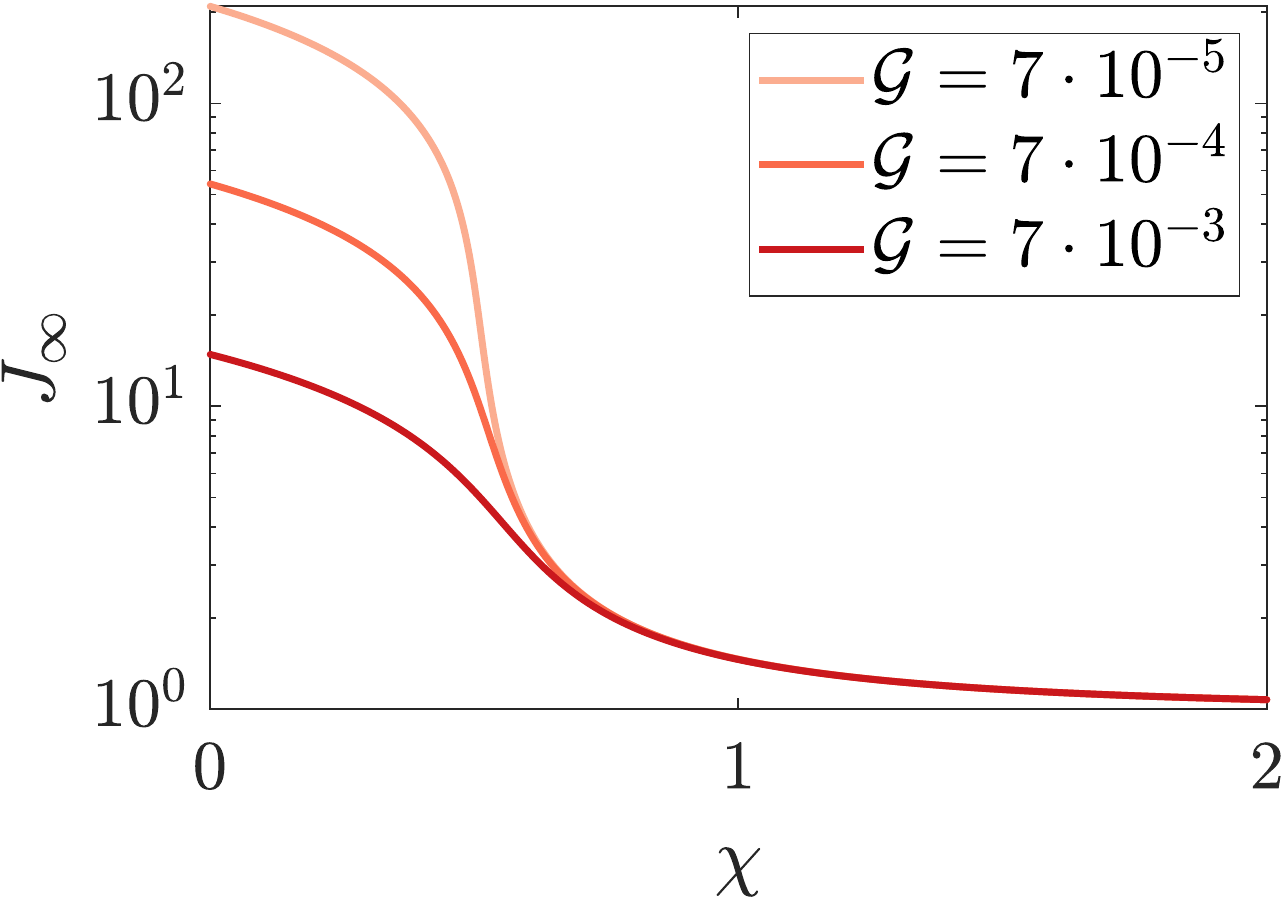}
  \includegraphics[width=0.32\textwidth]{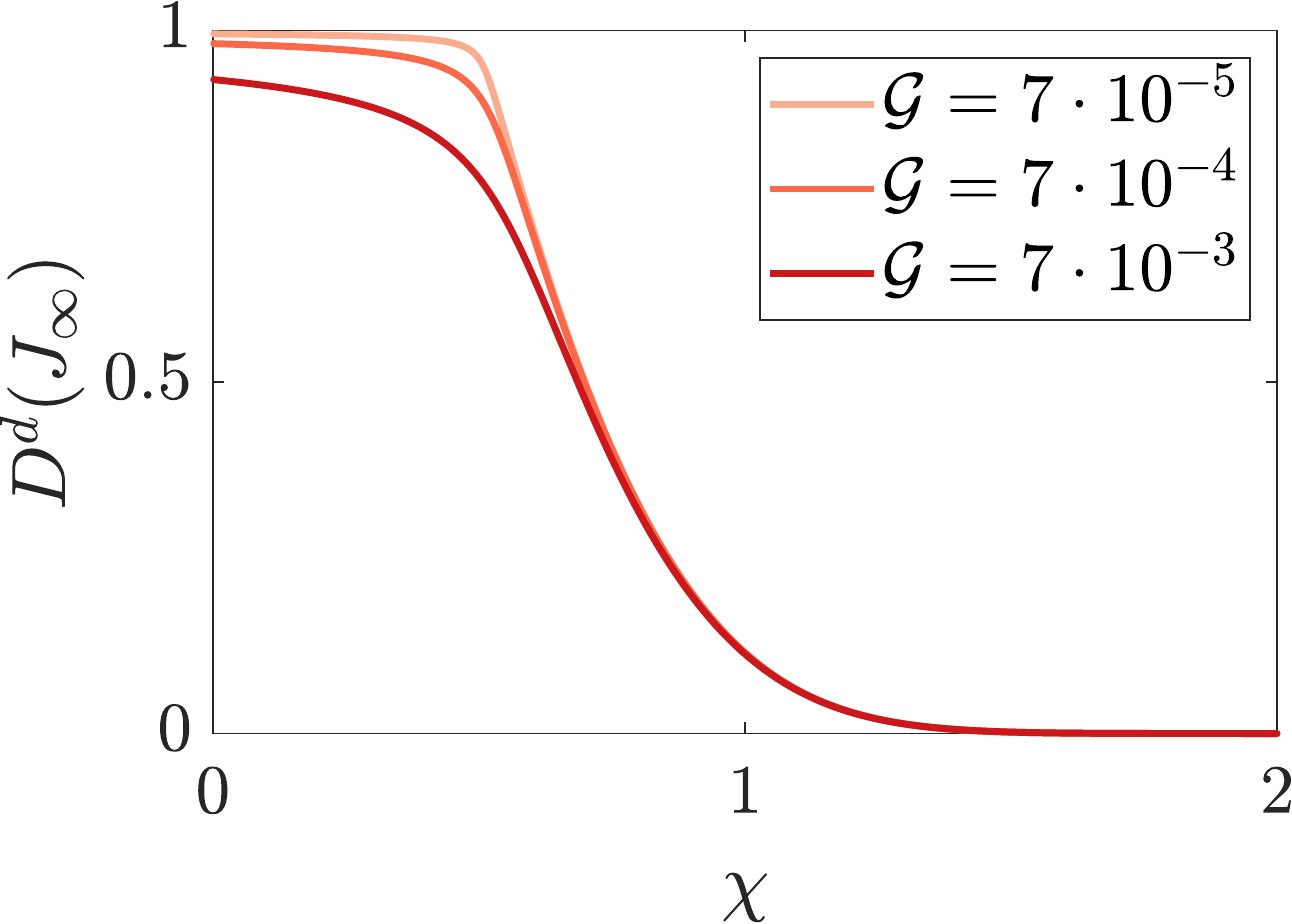}
  \includegraphics[width=0.32\textwidth]{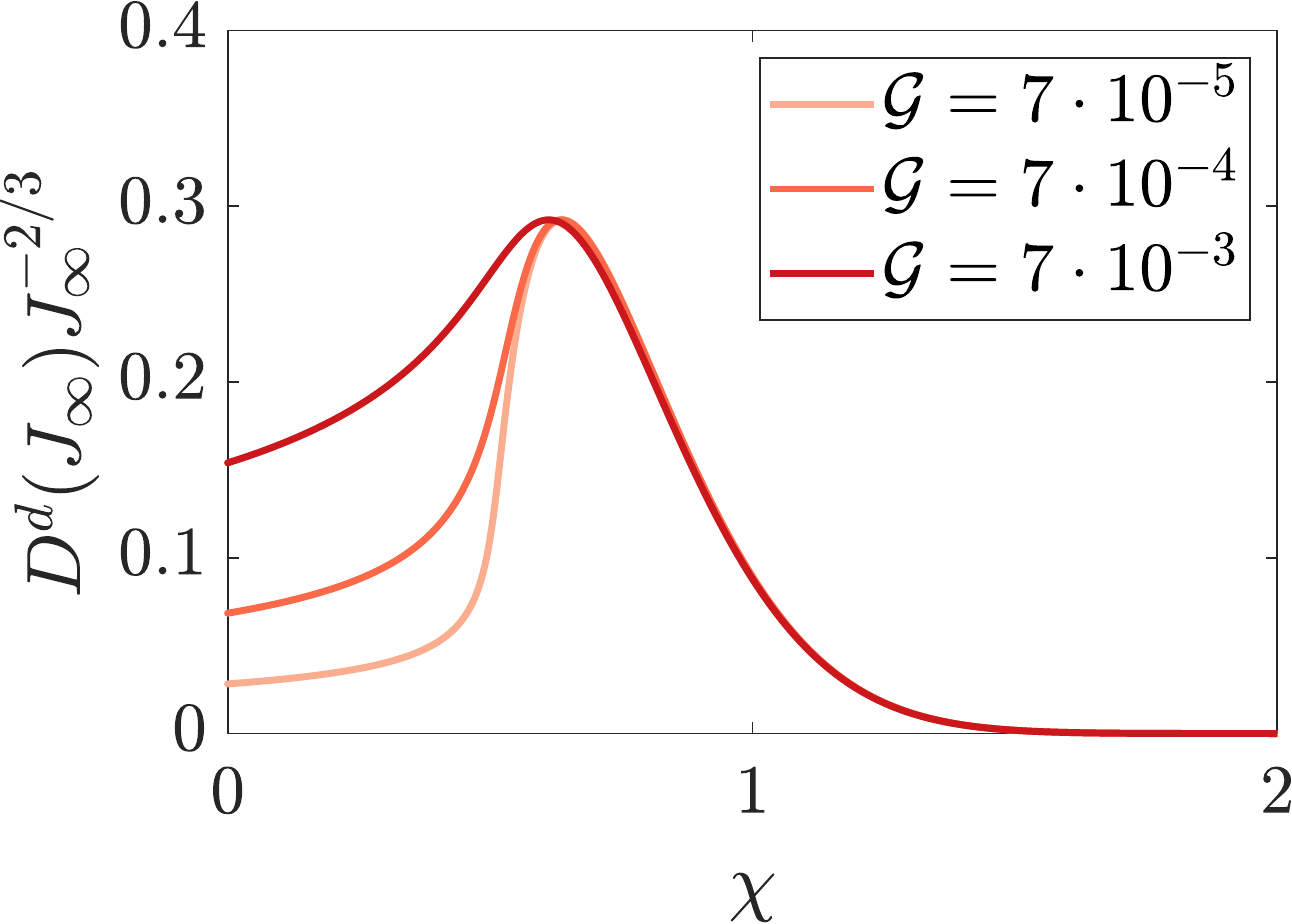}
  \caption{Equilibrium values of the (a) swelling ratio $J_\infty$,
    (b) drug diffusivity $\Dd(J_\infty)$, and (c)
    drug mobility $\Dd(J_\infty) J_\infty^{-2/3}$,
    plotted as functions of the Flory
    interaction parameter $\chi$ for different non-dimensional gel
    stiffnesses $\G$. The diffusivity $\Dd$ is given by \eqref{eqn:Dd}.}
  \label{fig:spherical_equilibria}
\end{figure}

As the Flory interaction parameter $\chi$ increases, there is a marked
decrease in the degree of swelling that occurs for all gel stiffnesses;
see Fig.~\ref{fig:spherical_equilibria}~(a). Consequently,
the equilibrium
drug diffusivity decreases as well;
see Fig.~\ref{fig:spherical_equilibria}~(b). 
The curves for the
equilibrium swelling ratio and drug diffusivity converge when $\chi > 0.7$,
indicating that elasticity no
longer plays a role in determining the equilibrium. For $\chi < 0.7$,
the degree of swelling becomes strongly dictated by the gel stiffness 
$\G$, with softer gels (smaller $\G$) undergoing larger deformations.
Correspondingly, the equilibrium drug diffusivity increases with decreasing
$\G$ as well.
For the softest gels ($\G = 7 \cdot 10^{-5}$), the drug diffusivity
$\Dd$ approaches unity, implying that drug
diffusion becomes uninhibited by the presence of the polymer network.

The increase in drug diffusivity due swelling is
counteracted by the increase in distance that drug molecules must travel
to reach the free surface. This increase in distance is captured by
the factor of $\lambda_r^{-2}$  in \eqref{nd:drug:eqn},
which is equal to $J_\infty^{-2/3}$
at equilibrium. We can thus define the
equilibrium drug mobility as $\Dd(J_\infty) J_\infty^{-2/3}$ to capture
the competing effects of increases in diffusivity and gel size. For
$0 < \chi < 0.7$, the equilibrium mobility increases with the gel
stiffness, as seen in Fig.~\ref{fig:spherical_equilibria}~(c). Thus,
stiffer gels are, in fact, more effective at releasing drug molecules
because the reduction in gel size completely offsets the smaller drug
diffusivity.
The strong dependence of the drug mobility on the gel size underscores the
importance of capturing finite deformations in the model. 

\begin{figure}
  \centering
  \subfigure[]{\includegraphics[width=0.48\textwidth]{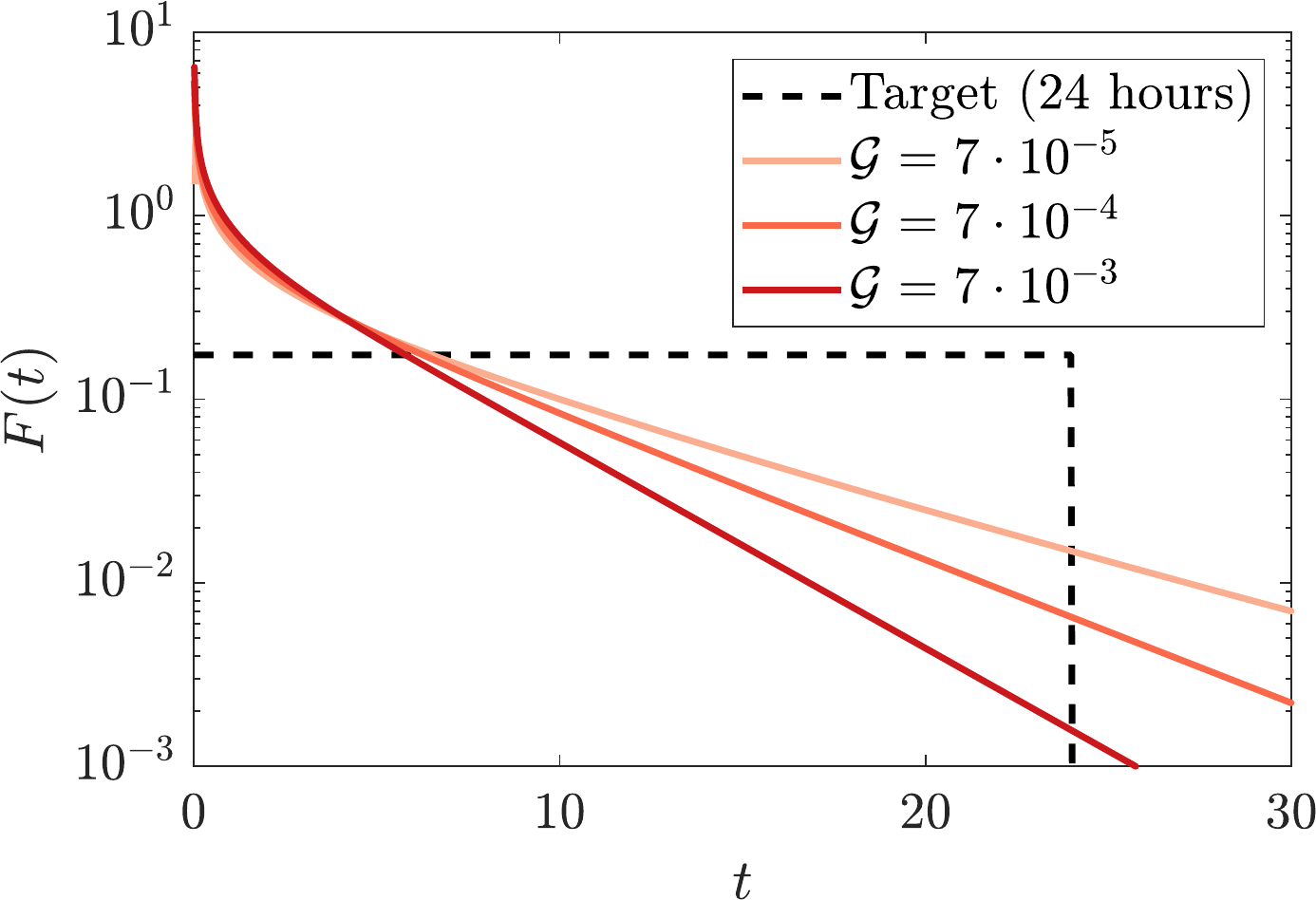}}
  \subfigure[]{\includegraphics[width=0.48\textwidth]{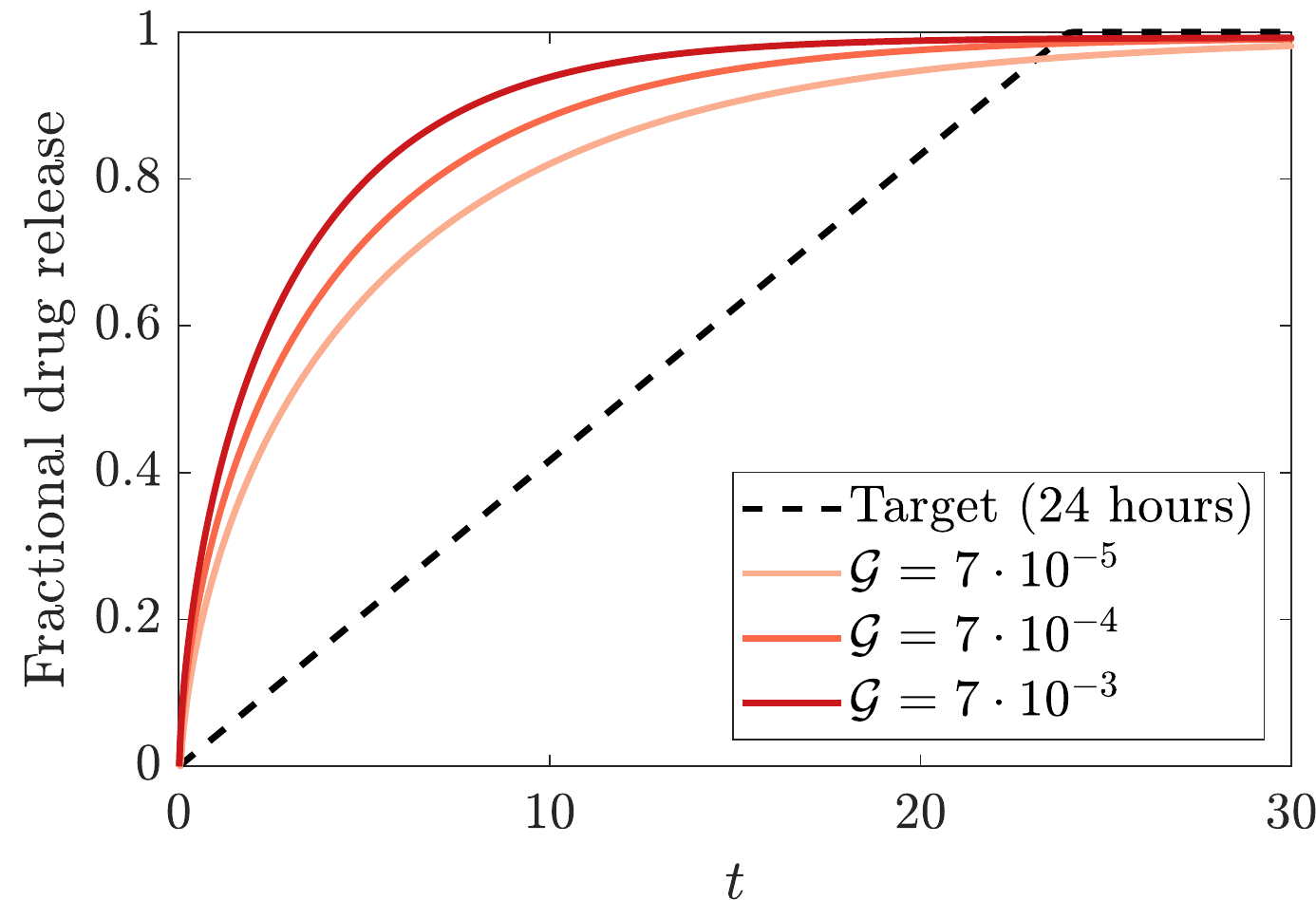}}
  \caption{Drug-release dynamics with a uniform initial loading of
    drug molecules,
    $d(R) \equiv 1$. (a) The efflux of drug
    molecules. (b) The fractional amount of drug released into
    the environment. The parameter values are $\chi = 0.5$ and
    $\D = 0.1$.}
  \label{fig:uniform}
\end{figure}

We explore the baseline drug-release dynamics for various gel stiffnesses
by assuming that the initial
loading of the drug is uniform; that is, we take $d(R) \equiv 1$.
The Flory parameter is set to $\chi = 0.5$ in order to capture changes in
the equilibrium state with gel stiffness. 
By numerically solving the governing equations, we compute the drug efflux
$F(t)$ according to \eqref{eqn:F} and the fractional drug release defined as
\begin{align}
  \mathcal{R}(t) = 1 - \frac{\mathcal{M}(t)}{\mathcal{M}(0)},
  \quad \mathcal{M}(t) = 4 \pi \int_{0}^{1} \Nd(R,t)R^2\, \d R.
  \label{eqn:R}
\end{align}
The numerical results are shown in Fig.~\ref{fig:uniform}, where, for
reference, we also plot a target profile in which the
drug efflux is constant over a 24-hour period, corresponding to a
piecewise-linear fractional drug release.
We consistently
observe a large initial ``burst'' where the drug efflux is large, resulting in
a rapid release of drug from the gel.
In particular, half of the drug molecules are released within
the first four hours.
The softest gels have the
slowest drug-release kinetics, as expected from the equilibrium drug mobility
seen in Fig.~\ref{fig:spherical_equilibria}~(c),
and thus provide the closest approximation to the target profile.

\section{Theory of optimal drug loading}
\label{sec:optimisation}

The bursting observed in Fig.~\ref{fig:uniform} is a generic phenomenon that occurs for all parameter values.  Achieving a uniform efflux, or a linear drug-release profile, simply by parameter variation is not possible.  Therefore, we explore how the initial drug concentration $d(R)$ can be varied to bring the
efflux $F(t)$ as close as possible to a prescribed target profile $A(t)$.

\subsection{Problem definition}

The objective function that measures the difference between $F(t)$ and $A(t)$ is taken to be
\begin{equation}
  \label{eq:OptimObj}
  H(d) := \int_0^{\infty} (F(t) - A(t))^2\, \d t,
\end{equation}
which penalises sustained deviations from the target profile; short periods of time where the drug release is different from
the desired profile are unimportant.
A similar objective function was used by Lu \etal\cite{lu1998modeling},
but here we ignore the cost of the drug, which is constant for a fixed value
of $\epsilon$, corresponding to a fixed amount of drug. 
While \eqref{eq:OptimObj} is a good measure of the ``closeness'' of $F(t)$ and $A(t)$, it is also important to impose that the total drug delivered is the same as the target total amount. This means
\begin{equation}
\label{eq:DosageConstraint}
    \int_0^{\infty} F(t)\, \d t = \int_0^{\infty} A(t)\, \d t = \frac{4}{3}\pi.
\end{equation}
The only other constraint on the function $d(R)$ is that
\begin{equation}
\label{eq:OptimConstr}
    d(R) \geq 0, \quad \forall R \in (0,1).
\end{equation}
The endpoint value $ d(1) = 0$, which is necessary for the boundary conditions, can be ignored as it contributes zero volume of drug to the system. For the moment, we will consider an infinitely dilute drug with
$\epsilon \to 0$ and thus will not impose an upper bound on
$d$.
Thus, the aim is to minimise (\ref{eq:OptimObj}) subject to (\ref{eq:DosageConstraint}) and (\ref{eq:OptimConstr}).

\subsection{Impossibility of a perfect solution}
\label{sec:impossible}
It may seem that there is enough freedom in choosing $d(R)$ to ensure that $F(t) \equiv A(t)$ for each possible choice of target profile $A$. However, the long-time behaviour of drug release means that this is not the case. Recall that, at long times, the hydrogel expands to a uniform equilibrium, as discussed in Sec.~\ref{sec:benchmarking}.
Thus, by assuming the convergence  to the equilibrium state is uniform,
then, to leading order, the diffusion equation for the drug becomes
\begin{equation}
\label{eq:constdiff}
    \frac{\partial N^d}{\partial t} \sim  \frac{\D \Dd(J_\infty) J_\infty^{-2/3}}{R^2}\frac{\partial}{\partial R}\left(R^2\frac{\partial N^d}{\partial R}\right) \text{    as } t \to \infty.
\end{equation}
By seeking a separable solution, it is straightforward to show that the leading behaviour in time of the efflux
is given by
\begin{equation}
F(t) =  -\td{}{t} \left(4 \pi \int_0^1 R^2 N^d(R,t)\,\d R\right) =  O\left(\exp \left(-\D \Dd(J_\infty) J_\infty^{-2/3} \pi t\right)\right) \text{  as } t \to \infty.
\label{eq:Decay}
\end{equation}
Thus, if the target profile $A(t)$ has slower-than-exponential decay, for
example, then it is impossible that $F(t) = A(t)$ everywhere, so any optimal solution will have a non-zero value of $H$.

\subsection{Formulation of a discrete optimsation problem}
The linearity of the drug-diffusion problem \eqref{nd:drug:all} can be
exploited to derive a discrete optimisation problem
that approximates the full problem given
by \eqref{eq:OptimObj}--\eqref{eq:OptimConstr}.
First suppose $\{N_i^d(R,t)\}_{i}$ is a set of $M$ functions satisfying \eqref{nd:drug:eqn} and its boundary conditions \eqref{nd:drug_bc}, with initial conditions $\{\xi_i(R)\}_i$. Moreover, let  $\{f_i(t)\}_i$ be the corresponding drug effluxes. Suppose further that $d_i$ are real constants. Then, by defining the initial drug concentration as
\subeq{
  \begin{align}
    d(R) &= \sum_{i=1}^M d_i \xi_i(R), \label{super:d}
  \end{align}
  the solution for the drug concentration $\Nd$ and the drug efflux $F$
  can be written as
  \begin{align}
    N^d(R,t) = \sum_{i=1}^M d_i N^d_i(R,t), 
    \quad 
    F(t) =  \sum_{i=1}^M d_i f_i(t). \label{super:F}
  \end{align}}
In light of \eqref{super:F}, we refer to each $f_i$ as a partial efflux.
The initial drug concentration
$d(R)$ in \eqref{super:d} is now written as a sum of
spherical delta functions that are centred at $M$ distinct radial coordinates
$0 < R_i < 1$
by defining
\begin{align}
  \xi_i(R) = \frac{\delta(R - R_i)}{4 \pi R^2}.
  \label{super:xi}
\end{align}
Each weight $d_i$ in
\eqref{super:d} therefore corresponds to the number of drug molecules located
at the point $R_i$. The discrete optimisation problem
will compute the optimal values for the $M$ weights
$d_1$, $d_2$, $\ldots$, $d_M$.
Since $f_i(t)$ denotes the drug efflux that is obtained from using $\xi_i$ as
an initial condition, we have that
\begin{align}
  \int_{0}^{\infty} f_i(t)\,\d t = 4 \pi \int_{0}^{1} \xi_i(R)R^2\, \d R = 1,
\end{align}
which implies that
\begin{equation}
  \int_0^{\infty}F(t)\,\d t = \sum_{i=1}^{M} d_i \int_0^{\infty}f_i(t)\, \d t = \sum_{i=1}^{M}d_i.
\end{equation}
By collecting the weights $d_i$ into a vector $\mathbf{d} = (d_1, d_2, \ldots, d_M)$, we can formulate
the following discrete optimisation problem:
\begin{equation}
    \label{eq:OverallOptim}
    \text{min}\left\{H(\textbf{d}) : \textbf{d} \geq \textbf{0} \text{  ,  } \sum_{i=1}^{M} d_i = \int_0^{\infty} A(t)\, \d t\right\},
\end{equation}
where the objective function $H$ is now
\begin{equation}
  H(\textbf{d})  = \int_0^{\infty}\left(\sum_{i=1}^{M}d_if_i(t)-A(t)\right)^2\,\d t, \label{dis:H}
\end{equation}
which is a quadratic function of the variables $d_i$. To practically calculate the integral in \eqref{dis:H}, it is necessary to restrict the domain of integration to $[0,T]$ for some large $T$; we find that $T = 50$ is a sensible choice.

\subsection{Convexity of the objective function}
\label{sec:convexity}

A useful property of the objective function $H$ given by \eqref{dis:H} is that it is a convex, quadratic function of $\textbf{d}$. This can be seen by expanding the integrand to give
\begin{equation}
\label{eq:HasQuad}
    H(\textbf{d}) = \sum_{i=1}^{M}\sum_{j=1}^{M}\left(d_id_j\int_0^{\infty}f_i(t)f_j(t)\, \d t\right) -2\sum_{i=1}^{M}\left(d_i\int_0^{\infty} A(t)f_i(t)\, \d t\right) + \int_0^{\infty} A(t)^2\, \d t,
\end{equation}
which can be concisely written as
\begin{equation}
    H(\textbf{d}) = \frac{1}{2}\textbf{d}^T \boldsymbol{S} \textbf{d} - \textbf{d}^T\textbf{q} +\int_0^{\infty} A(t)^2\, \d t ,
\end{equation}
where $\boldsymbol{S} \in \R^{M \times M}$ and $\textbf{q} \in \R^{M}$ are defined by
\begin{equation}
    S_{ij} = 2\int_0^{\infty}f_i(t)f_j(t) \,\d t, \qquad q_i = 2\int_0^{\infty} A(t)f_i(t)\, \d t.
\end{equation}
In particular, $\boldsymbol{S}$ is positive semi-definite as, for any $\mathbf{v}\in\R^{M}$,
\begin{align}
  \mathbf{v}^T \boldsymbol{S} \mathbf{v}
    &= 2\int_0^{\infty} \sum_{i=1}^M\sum_{j=1}^M\left(v_if_i(t)f_j(t)v_j\right)\, \d t
    = 2\int_0^{\infty} \left(\sum_{i=1}^Mv_if_i(t)\right)^2\, \d t \geq 0.
\end{align}
Therefore, $H$ is convex as $\boldsymbol{S}$ is the Hessian matrix of $H$.
Using standard convex programming results~\cite{niculescu2006convex}, we can prove that a constrained global minimum exists and that
constrained local and global minima are equivalent; details
of the proofs are provided in Appendix~\ref{app:optim}.
Thus, the global minimum can be calculated by simply finding
a local minimum.  For the remainder of this paper, the optimal value
of $H$ will be denoted as $H^*$.

\subsection{Numerical implementation}

The discrete optimisation problem \eqref{eq:OverallOptim}--\eqref{dis:H} is numerically solved by discretising the governing equations using the finite-difference method described in Sec.~\ref{sec:numerics}.  In particular, the radial domain is discretised into $M$ cells of width $\Delta R = 1 / M$.  The radial coordinates $R_i$ used to formulate the discrete optimisation problem (see \eqref{super:xi}) are chosen to coincide with the positions of the $M$ cell midpoints. The domain of each cell can then be defined
as $\Omega_i = \{R: R_{i}^{-} < R < R_{i}^{+}\}$ where $R_{i}^{\pm} = R_i \pm \Delta R / 2$ represent the cell edges. 
The spherical delta functions given by \eqref{super:xi} are replaced with step functions defined by
\begin{align}
  \xi_i(R) = \begin{cases}
    \frac{3}{4 \pi}\left[(R_{i}^{+})^3 - (R_{i}^{-})^3\right]^{-1}, &\quad R \in \Omega_i, \\
    0, &\quad \text{otherwise},
  \end{cases}
         \label{num:xi}
\end{align}
where  $i = 1, 2, \ldots, M$. The functions $\xi_i$ in \eqref{num:xi} can
be interpreted as localised packets of drug located at the $i$-th cell.

For a given set of parameter values,
we first numerically solve the dimensionless
hydrogel equations. We then numerically solve the drug-diffusion problem using each initial condition $\xi_i$ in \eqref{num:xi} to compute the drug concentration $\Nd_i$ and the partial efflux $f_i$. The set of partial effluxes $\{f_i\}_i$ is then used in the discrete optimisation problem
\eqref{eq:OverallOptim}--\eqref{dis:H},
which is solved using MATLAB's quadratic programming algorithm {\tt quadprog}.  The analytical results developed in Sec.~\ref{sec:convexity} and Appendix~\ref{app:optim} ensure that {\tt quadprog} will rapidly converge to a global minimum. The ability to
pre-compute the partial effluxes, which requires one solution to the nonlinear hydrogel model and $M$ solutions of the linear drug-diffusion model, results
in a highly efficient scheme for the numerical optimisation.
After determining the optimal weights $d_i$, the total efflux $F$ and the drug concentration $\Nd$ can be constructed using \eqref{super:F}. 

\subsection{The partial effluxes from localised drug loadings}

\begin{figure}
    \centering
    \includegraphics[width=0.48\textwidth]{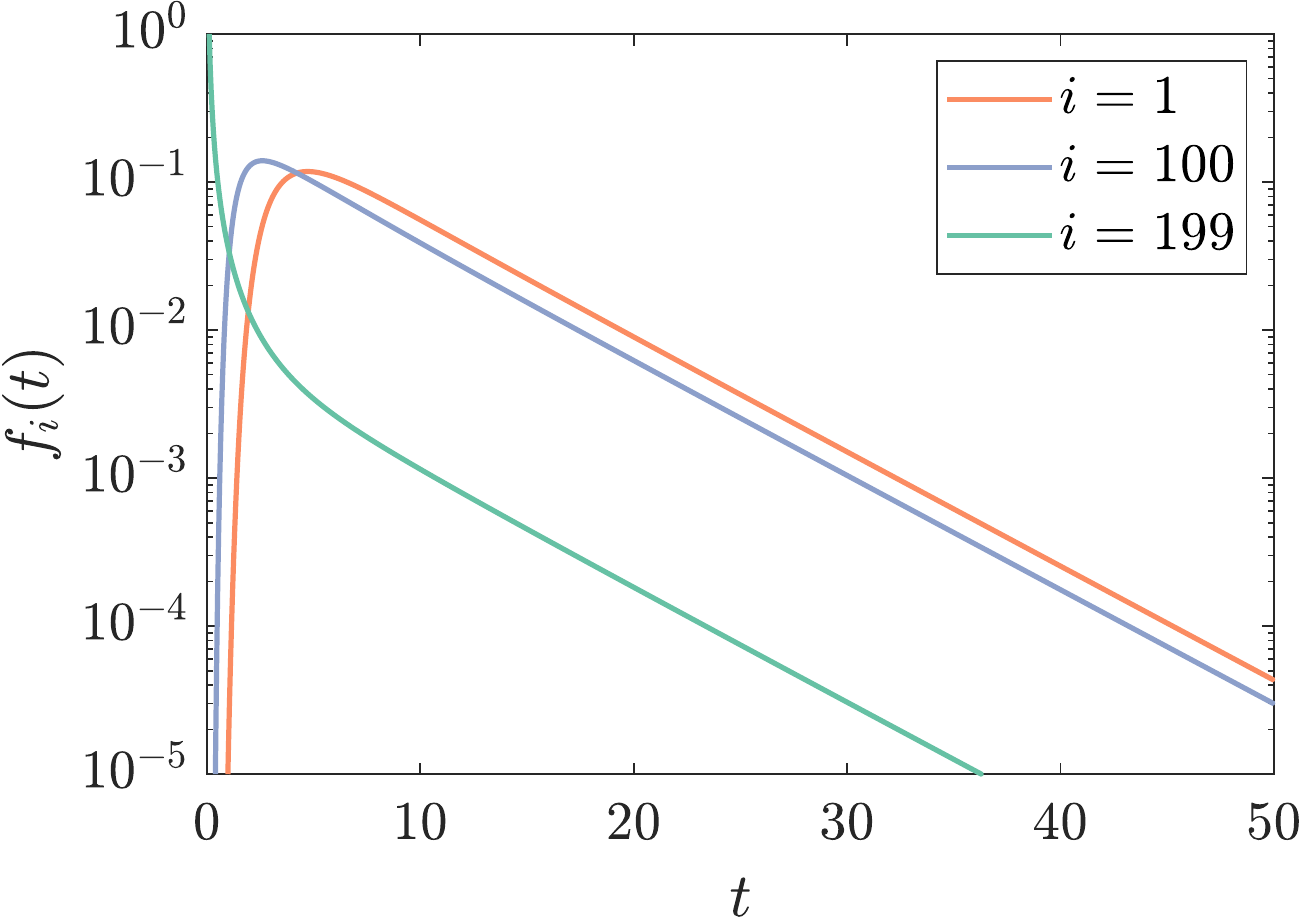}
    \caption{Plots of the partial drug effluxes $f_i(t)$ for $i=1$, $i=100$, and $i=199$ due to the packet-like initial concentrations \eqref{num:xi}.
      The values of $i$ correspond to
      $R=0$, $R=0.5$ and $R=1$, respectively.
      The parameter values are $\chi = 0.5$, $\D = 0.1$, and $\G = 7\cdot 10^{-4}$. The numerical simulations are based on using $M = 199$ cells.}
    \label{fig:fiForms}
\end{figure}

The partial effluxes $f_i(t)$ obtained from the localised drug packets given
by \eqref{num:xi}
act as basis functions in the construction of the total efflux $F(t)$, as seen
in \eqref{super:F}. 
Therefore, the key to understanding the optimal solution lies in the time evolution of the partial effluxes $f_i(t)$.
If a drug packet is placed sufficiently close to the
free boundary of the gel, then the correponding efflux $f_i$
monotonically decreases from a large initial value; see Fig.~\ref{fig:fiForms}.  In this case, the large initial efflux is driven by the incompatibility between the initial condition and the perfect sink boundary condition,
the latter of which forces the drug concentration to rapidly tend to zero.
If a drug packet is placed in the bulk of the gel, then the efflux first increases from zero to a peak value, after which it exponentially decays.  The transient increase in efflux is driven by the increase in drug mobility that occurs due to swelling.  As placement of the drugs moves further from the outer boundary,  the longer it takes for the efflux to reach its peak value,
as can be seen from Fig.~\ref{fig:fiForms}.  Eventually, all of the effluxes $f_i$ exponentially
decay with the same rate, highlighting the unavoidable long-term behaviour of the drug efflux discussed in Sec.~\ref{sec:impossible}.  

From the observations made from Fig.~\ref{fig:fiForms}, it follows that
the drug concentration close to the boundary can be utilised to capture the target profile $A(t)$
at small times.  The drug concentration in the bulk of the gel enables the target
profile to be captured at intermediate times.  Capturing the target profile after all of the effluxes
$f_i$ have peaked is particularly difficult and requires amplifying the drug concentration
near the gel centre at the expense of potential overshoots at intermediate times.

\section{Case studies}
\label{sec:case_studies}

We apply our theory of optimal drug loading to specific scenarios by considering
non-dimensional piecewise-constant target profiles with the same form as
\eqref{nd:A}.

\subsection{Optimal solutions for a 12-hour drug-release period}
\label{sec:12}

As a generic example, consider the optimal solution when $\tau = 12$,
corresponding to a constant release of drug over
a 12-hour period. The parameter values are set to  $\chi = 0.5$,
$\G = 7 \cdot 10^{-4}$, and $\mathcal{D} = 0.1$, which
are the same as those used when computing the partial effluxes $f_i(t)$ in
Fig.~\ref{fig:fiForms}.

\begin{figure}
  \centering
  \subfigure[]{\includegraphics[width=0.48\textwidth]{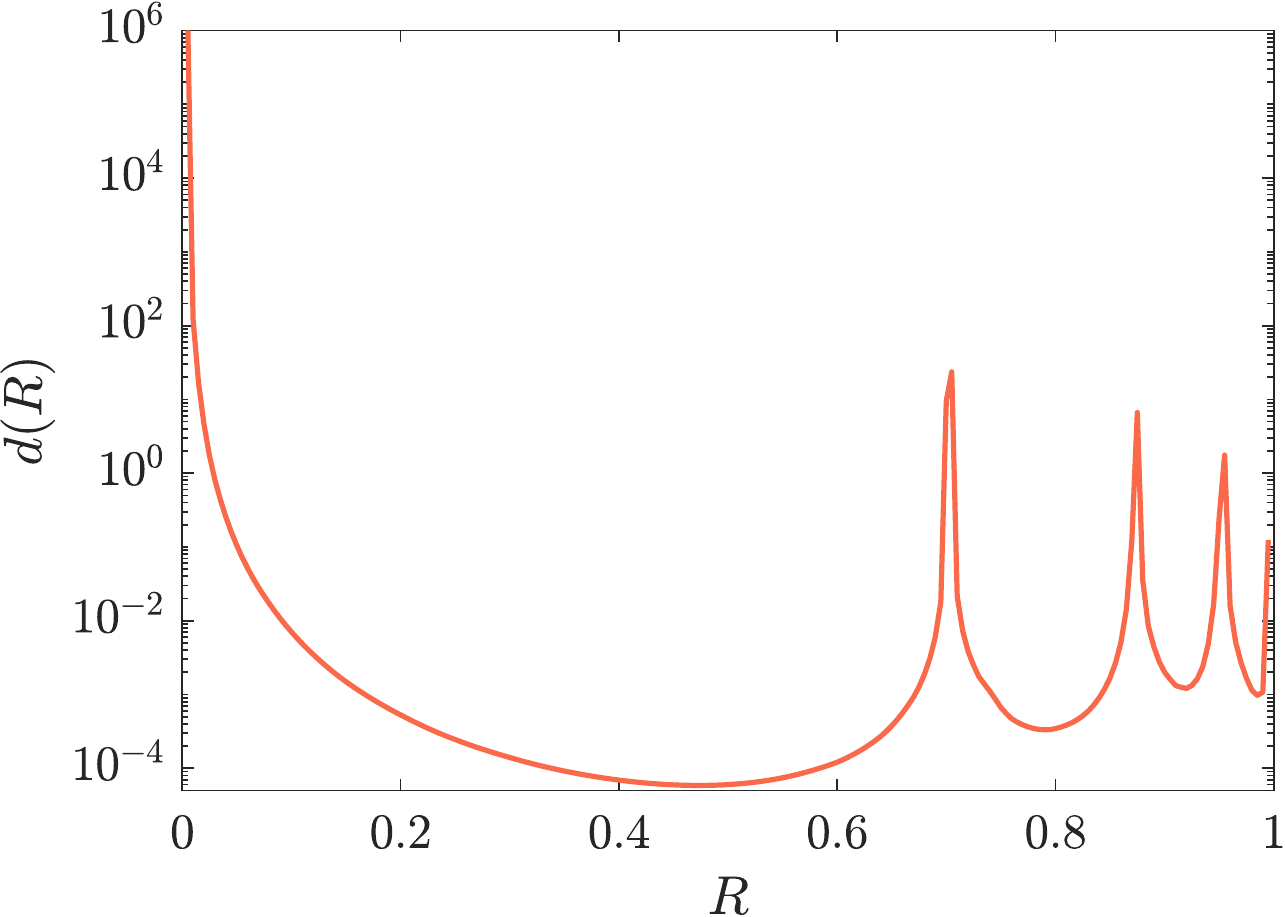}}
  \subfigure[]{\includegraphics[width=0.48\textwidth]{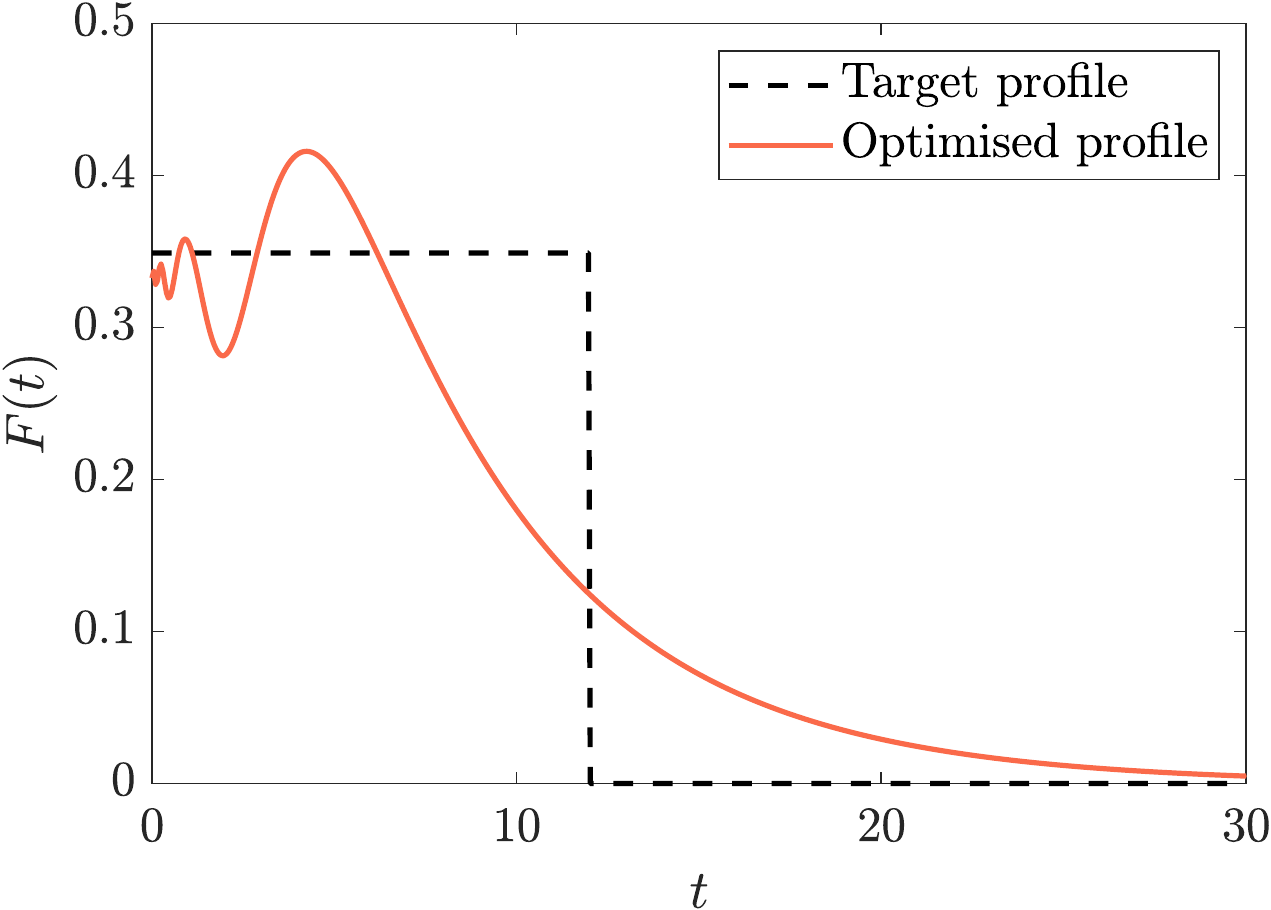}} 
  \caption{The optimal (a) initial drug concentration $d(R)$ and (b)
    corresponding drug efflux $F(t)$. The dashed line represents the
    piecewise-constant target profile. The parameter values are
    $\chi = 0.5$, $\D = 0.1$, $\G = 7 \cdot 10^{-4}$, and $\tau = 12$.}
  \label{fig:GeneralCaseDistrib}
\end{figure}

The optimal initial drug distribution is formed of distinct, concentrated packets that are separated by large drug-free regions; see Fig.~\ref{fig:GeneralCaseDistrib}~(a). The concentration in the
packets  decreases as their distance from the gel centre increases.
Such a distribution is expected on physical grounds.
Small concentrations of drug near the free boundary offset the initial
largeness of the corresponding partial effluxes seen in
Fig.~\ref{fig:fiForms}. The drug in the central packet at $R = 0$
sustains the long-term drug efflux. However, to reach the free surface,
the drug molecules in the central packet much diffusively spread across the
entire gel, resulting in a diminished concentration gradient and hence
diffusive flux. The largeness of the concentration in the central packet
offsets this behaviour.

The corresponding optimal drug efflux $F(t)$ consists of a sequence of
pulses of increasing amplitude; see Fig.~\ref{fig:GeneralCaseDistrib}~(b).
Each pulse is associated with one of the packets in the initial drug concentration.
The small quantity of drug in the packet near $R=1$ is responsible for the
first pulse, which provides an approximation of the target profile $A(t)$ for small times. Then, as the pulse from this packet of drug diminishes,
the pulse from the next packet begins, counteracting this decrease.
This pattern continues until the drug in the central packet nearest $R=0$
is released to create the largest pulse, which then exponentially decays. 
Due to the exponential tail,
the optimal efflux overshoots the target profile after the discontinuity
at $t = 12$. The overshoot
is compensated by a substantial undershoot beforehand when
$7 < t < 12$. The undershoot is itself compensated by an overshoot when
$4 < t < 7$, and the sequence repeats until $t = 0$.

\begin{figure}
\subfigure[]{\includegraphics[width=0.48\textwidth]{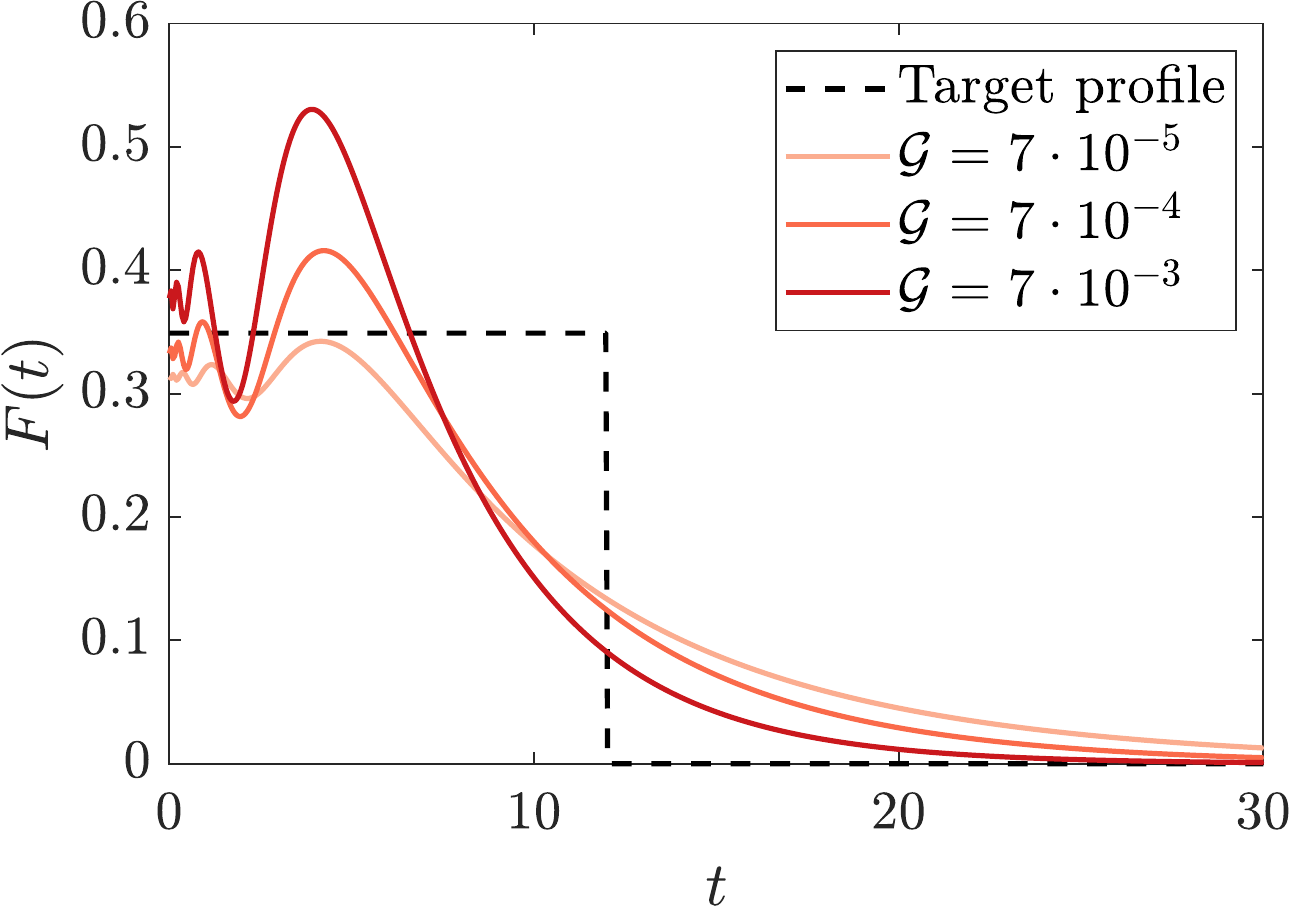}} 
\subfigure[]{\includegraphics[width=0.48\textwidth]{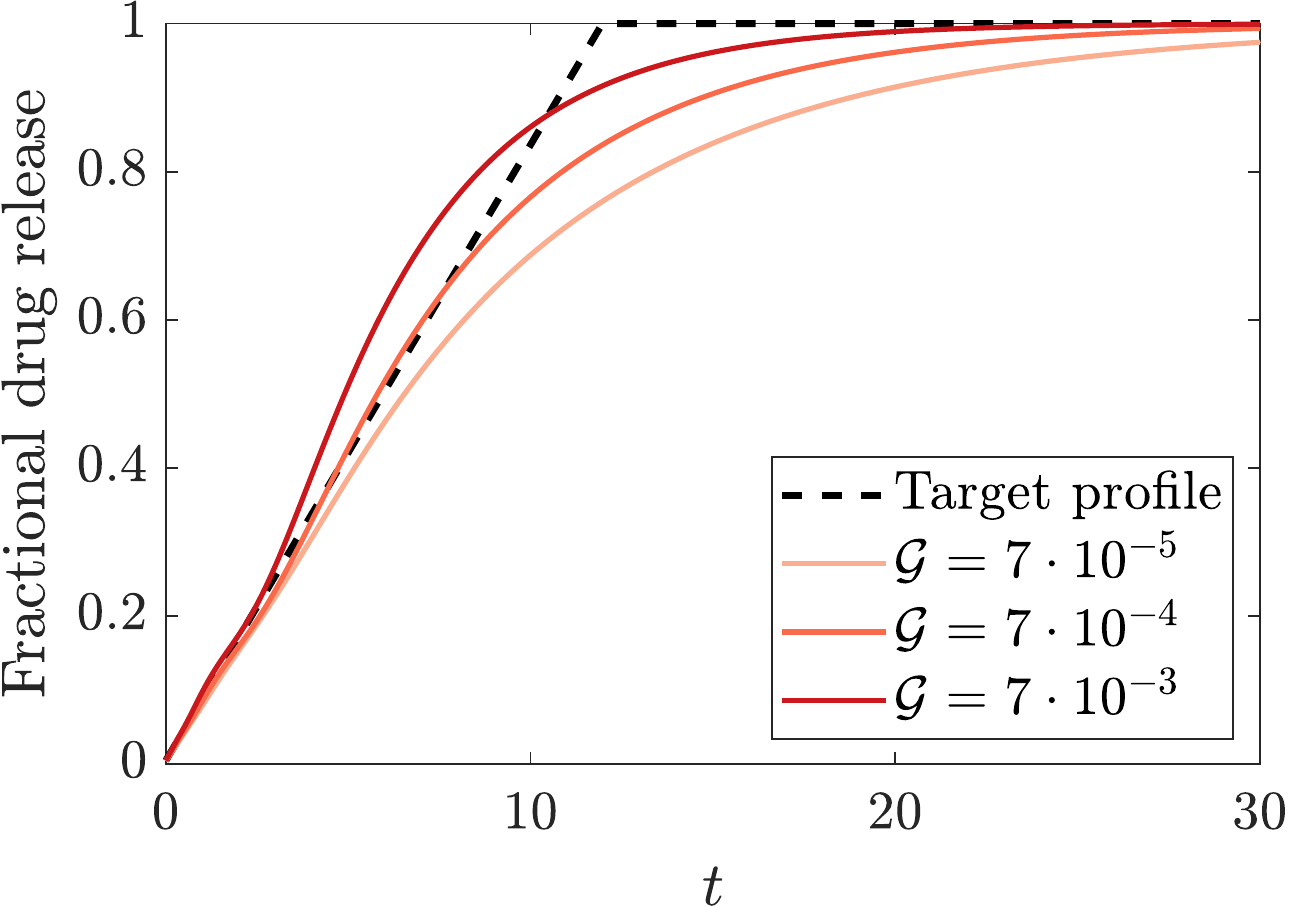}}
\caption{A comparison of optimal and target (a) effluxes and (b) fractional drug-release profiles for different gel stiffnesses $\G$ in the case of
  $\chi = 0.5$, $\D = 0.1$, and $\tau = 12$.}
\label{fig:GeneralCaseFluxes}
\end{figure}

The gel stiffness $\G$ plays an important role in the optimal solution
by modulating the equilibrium drug mobility 
$\Dd(J_\infty) J_\infty^{-2/3}$ and hence the decay rate of the efflux $F(t)$. 
When $\chi = 0.5$, the mobility increases
with the gel stiffness; see Fig.~\ref{fig:spherical_equilibria}~(c). 
For stiff gels with $\G = 7 \cdot 10^{-3}$,
the faster decay rate leads to less overshoot after the discontinuity
in the target profile, but results in much larger undershoots and overshoots
beforehand; see Fig.~\ref{fig:GeneralCaseFluxes}~(a). For soft gels with
$G = 7 \cdot 10^{-5}$, the slower decay rate leads to a greater
overshoot after the discontinuity, which is compensated by an efflux that
is consistently below the target profile beforehand.

By computing the fractional drug release using \eqref{eqn:R}, we find that the
gel stiffness can be combined with the optimal loading to tune the
drug-release profile and, in particular, eliminate the burst effect.
For an intermediate stiffness
of $\G = 7 \cdot 10^{-4}$, the target drug-release profile is perfectly captured
during the first 7.5 hours, as shown in Fig.~\ref{fig:GeneralCaseFluxes}~(b),
despite the sequence of undershoots and overshoots
that occur in the efflux. However, after 7.5 hours, there is a sharp
decrease in the release rate, resulting in a marked departure from the target
profile and a prolongation of the drug-release period. In particular, 
19 hours are needed for 95\% of the drug to be released.
By increasing the gel stiffness to $\G = 7\cdot 10^{-3}$,
the release rate can be accelerated,
which leads to a temporary overshoot compared to the target release profile but
lessens the long-term undershoot; in this case, only 14 hours are
required for 95\% of the drug to be released. Decreasing the
stiffness to $\G = 7 \cdot 10^{-5}$ leads to a slower drug-release profile that
is consistently below the target curve.



\subsection{Tuning the drug-release profile}
\label{sec:tuning}

Motivated by the results in Fig.~\ref{fig:GeneralCaseFluxes}, we now
explore how the gel stiffness can be used to further optimise the
drug-release profile. More specifically, we compute the optimal
value of the objective function $H$, denoted by $H^*$, across a range
of gel stiffnesses $\G$ and different combinations of the drug
diffusivity $\D$ and drug-release period $\tau$.

\begin{figure}
    \centering
    \subfigure[]{\includegraphics[width=0.32\textwidth]{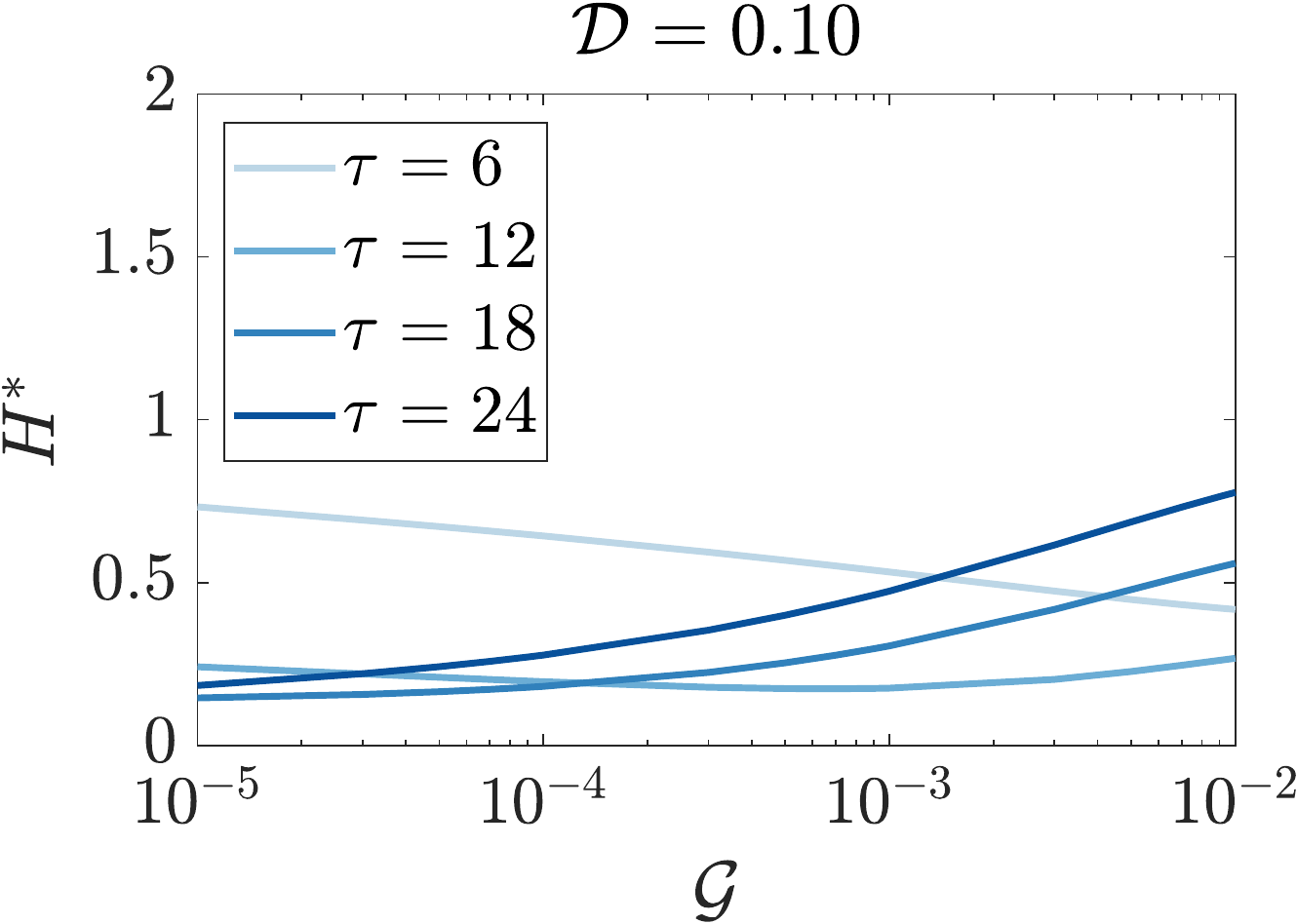}}
    \subfigure[]{\includegraphics[width=0.32\textwidth]{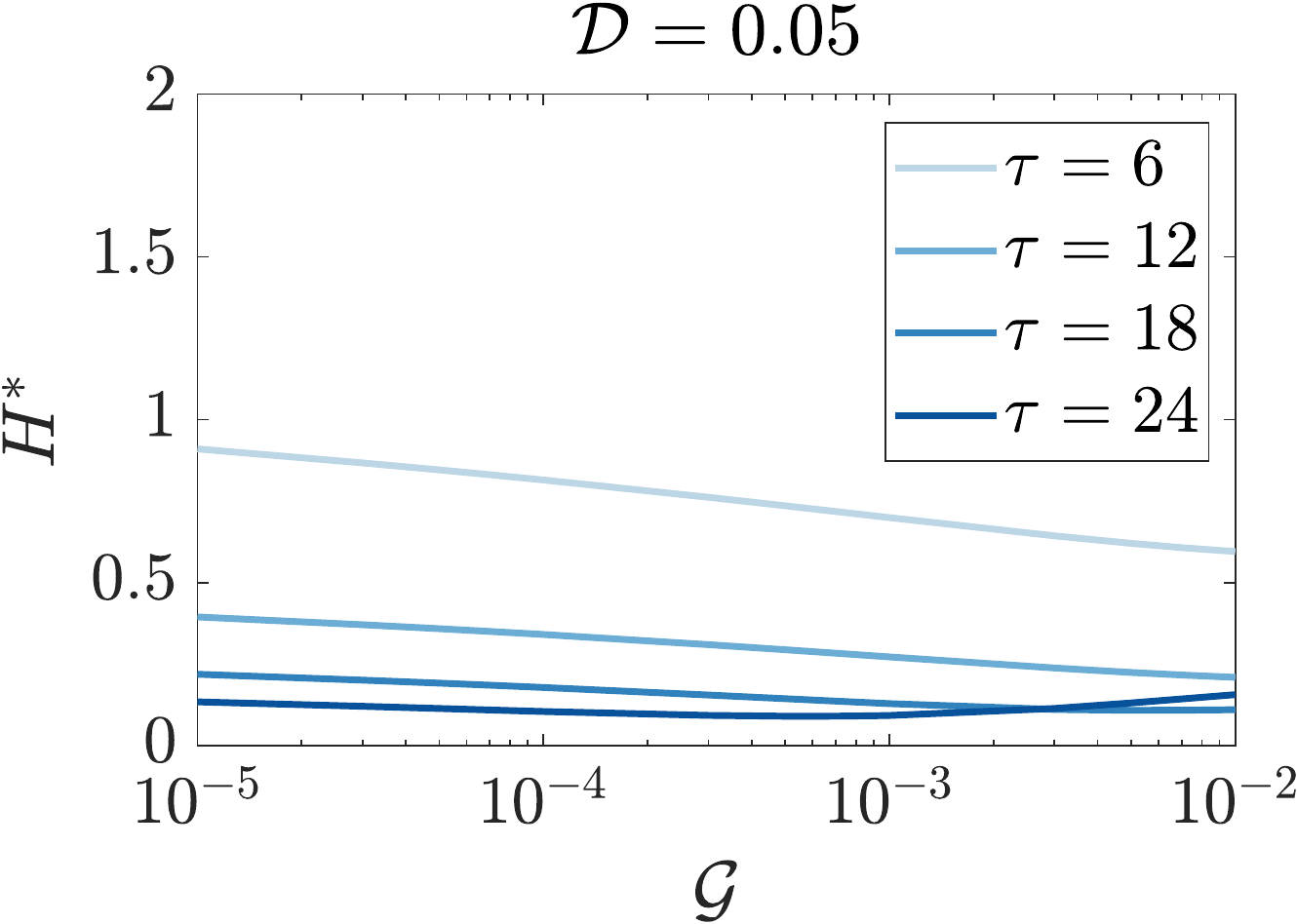}}
    \subfigure[]{\includegraphics[width=0.32\textwidth]{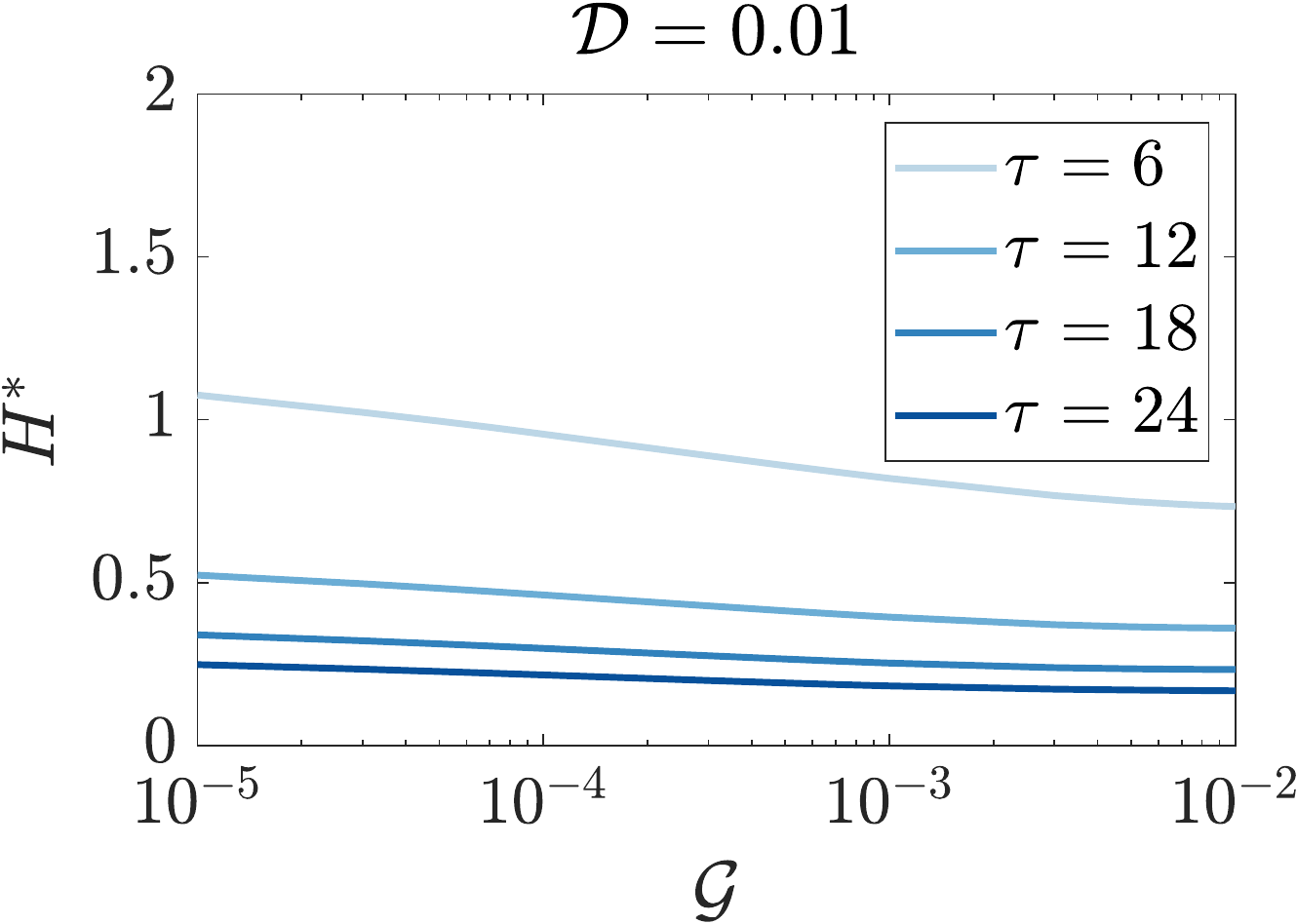}}
    \caption{The variation in the optimal value of the objective function
      $H^*$ with the gel stiffness $\G$ for different drug-release
      periods $\tau$ and non-dimensional drug diffusivities $\D$.
      Here, we take $\chi = 0.5$.}
    \label{fig:HG}
\end{figure}

We first consider the case when $\D = 0.1$. For a 12-hour drug-release
period with $\tau = 12$, the curve of $H^*$ as a function of $\G$ has
a global, internal minimum at $\G = 7 \cdot 10^{-4}$; see
Fig.~\ref{fig:HG}~(a). Thus, the target profile is best
approximated when the gel stiffness is $7 \cdot 10^{-4}$, in agreement with
the results in Fig.~\ref{fig:GeneralCaseFluxes}. Hydrogels that
are stiffer or softer would increase or decrease the equilibrium drug
mobility, respectively, and hence lead to drug molecules that are released
too quickly or slowly to be optimal. The increase in drug-release rate
that occurs for stiffer hydrogels can be advantageous for capturing
target profiles with smaller drug-release periods.
Conversely, softer hydrogels, with their slower drug-release kinetics,
will be better suited for capturing target profiles with larger drug-release
periods. Indeed, when the
drug-release period $\tau$ is decreased to 6 hours, the objective
function $H^*$ monotonically decreases with $\G$; when $\tau$
is increased to 18 or 24 hours, the
objective function monotonically increases with $\G$.

Decreasing the drug diffusivity via the dimensionless parameter $\D$
leads to stiffer gels performing better, as seen in
Fig.~\ref{fig:HG}~(b)--(c). In this case, the reduction in swelling and the
smaller radial extent of the hydrogel makes up for the decrease in the
rate of drug diffusion. 
As a final remark, the model can be simplified
for small values of $\D$ by rescaling time as $t = \D^{-1} t'$ and then
taking the limit $\D \to 0$. The resulting quasi-static model describes
linear drug diffusion throughout a uniformly swollen hydrogel that is in
chemo-mechanical equilibrium with the surrounding environment.

\subsection{Optimal solution with finite drug diluteness}

The numerical results in Sec.~\ref{sec:12}--\ref{sec:tuning}
are based on the assumption
of an infinitely dilute drug, as characterised by the limit $\epsilon \to 0$.
However, in practice, the drug molecules can account for roughly
10\% ($\epsilon = 0.1$) of the total initial volume~\cite{caccavo2015}.
Thus, when computing the optimum initial volume fraction of drug,
$\phid(R,0) = \epsilon d(R)$, using the solution for $d$ shown
in Fig.~\ref{fig:GeneralCaseDistrib}~(a), we see that it exceeds one.
This unphysical result stems from not imposing upper bounds on the
initial drug concentration when working in the infinitely dilute limit. 


We now consider situations where the drug molecules have
a finite diluteness. Thus, the global drug fraction $\epsilon$ is taken
to be a small but finite number.
We assume that the initial volume fraction of drug fraction satisfies
$0 < \phid(R,0) < 0.2$, which implies that $0 < d(R) < 0.2 \epsilon^{-1}$.
The optimisation theory developed in Sec.~\ref{sec:optimisation} still
applies because the
upper bound on $d(R)$ translates into upper bounds on each $d_i$ in the
discrete optimsation problem, which are straightforward to accommodate; see
Appendix~\ref{app:optim} for full details. 
Thus, the newly constrained
discrete optimisation problem admits global minima that can be readily computed
using Matlab's {\tt quadprog} function.

\begin{figure}
  \centering
  \subfigure[]{\includegraphics[width=0.48\textwidth]{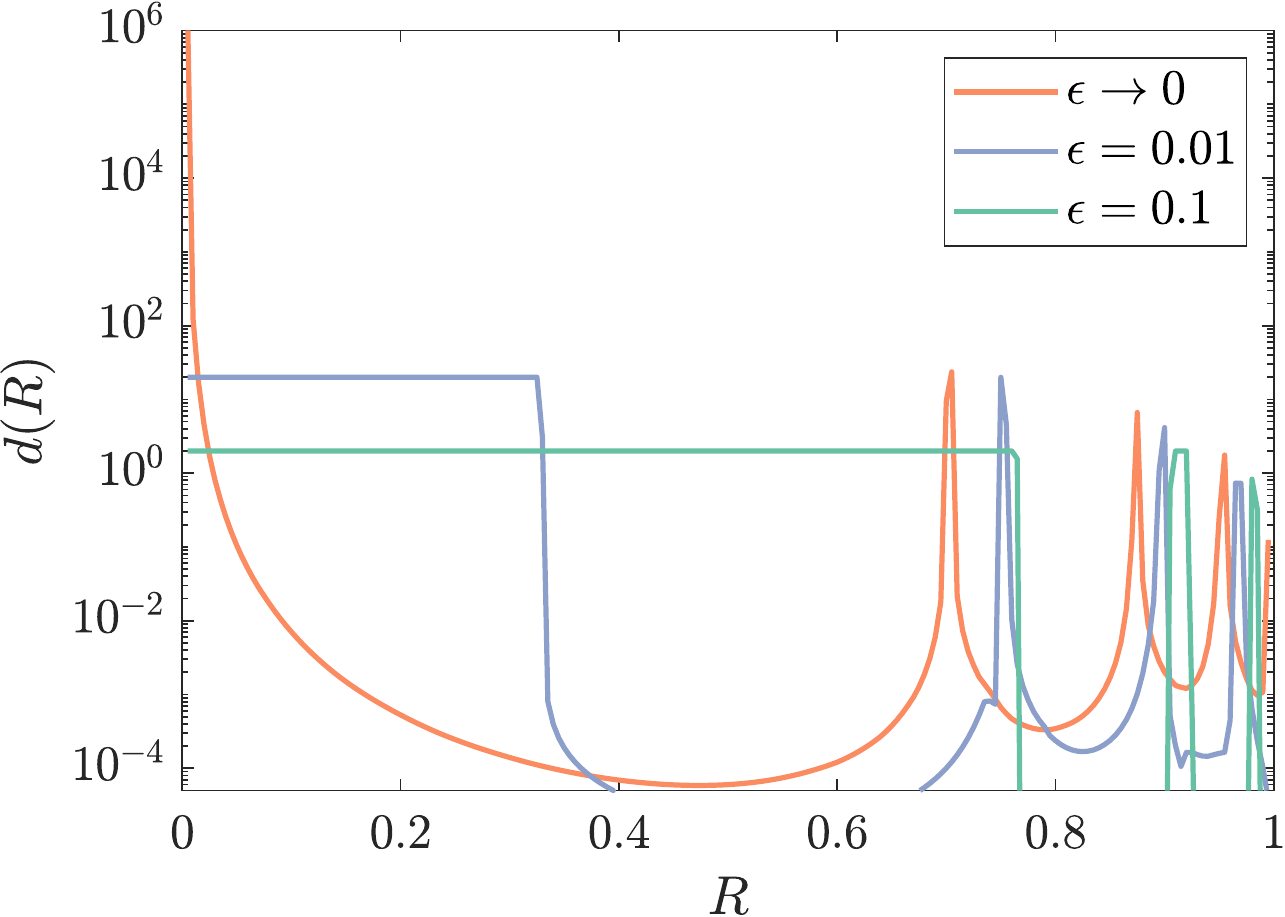}}
  \subfigure[]{\includegraphics[width=0.48\textwidth]{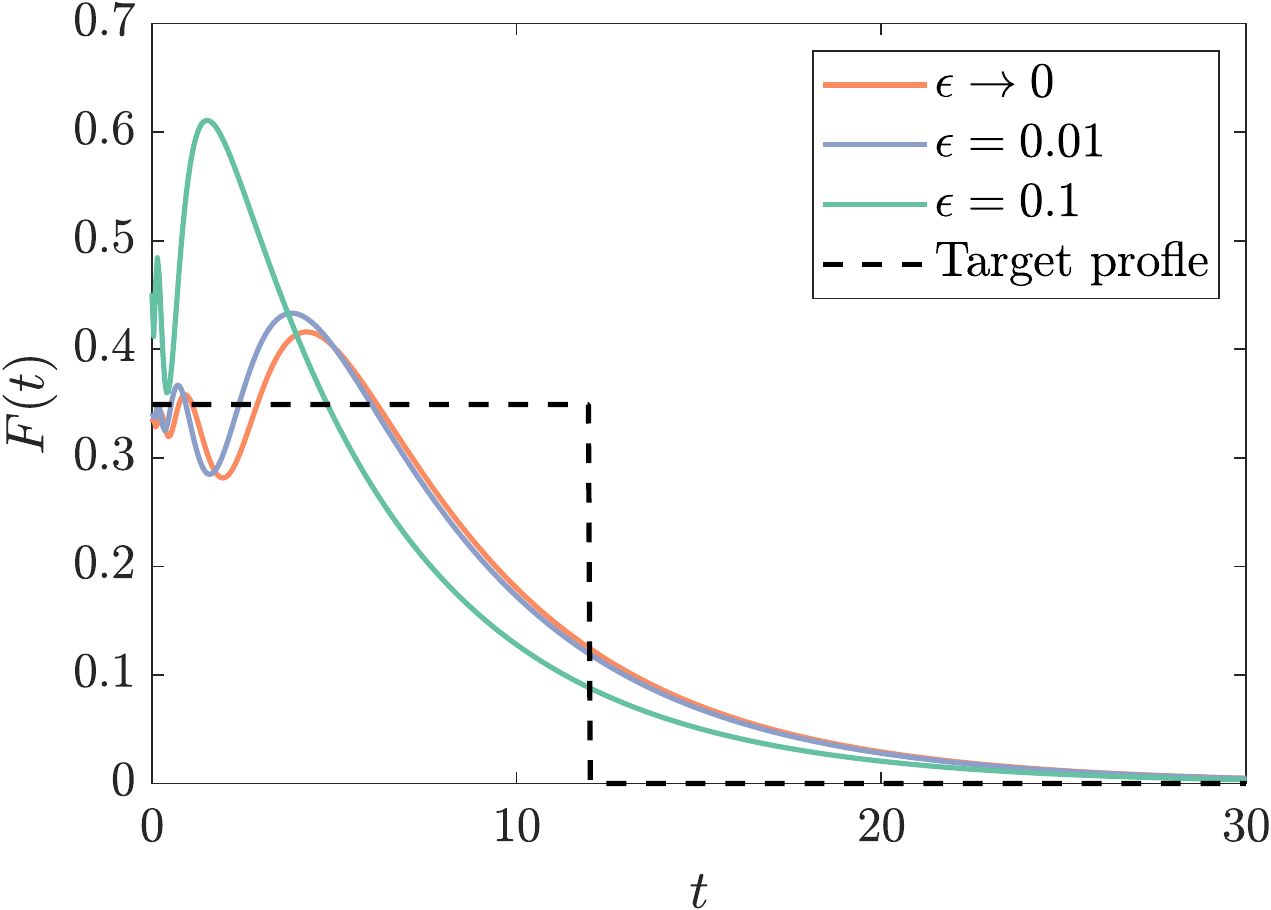}}
  \caption{Optimal (a) initial drug concentrations $d(R)$ and (b) effluxes
    $F(t)$ for different maximal drug dilutions, as measured through the
    global drug fraction $\epsilon$ defined by \eqref{eqn:d_cons}. The
    parameter values are $\chi = 0.5$, $\G = 7 \cdot 10^{-4}$, $\D = 0.1$.}
  \label{fig:eps_dF}
\end{figure}

As the global drug fraction $\epsilon$ is increased from zero,
the optimal initial drug concentration $d(R)$ is found to retain a
structure that consists of several drug-loaded packets near the free boundary
of the gel;
see Fig.~\ref{fig:eps_dF}~(a). 
However, the central packet near $R = 0$
that was
observed in the infinitely dilute case ($\epsilon \to 0$) has been replaced
with a uniformly loaded core in which the local drug fraction takes on
its maximum value. The radial extent of the drug-loaded core increases with
the global drug fraction $\epsilon$. 
To explain these results, we recall that the
initial drug concentration becomes increasingly constrained as $\epsilon$
increases from zero. Thus, if the initial drug concentration for the
infinitely dilute case has any packets that exceed the maximum allowable
concentration, then the drug molecules in these packets
are simply distributed over a greater
volume, that is, a over greater radial extent.
Any packets that have concentrations below the maximum are mostly
unaffected by the constraint, although their position might shift slightly.
From the optimal initial drug concentrations, we can conclude that
when increasing the total drug load for a fixed drug-release period $\tau$,
it is advantageous to preferentially
place drug molecules at the centre of the hydrogel.

The similarities in the optimal initial drug profiles
lead to optimal drug effluxes $F(t)$ that still consist of pulse
sequences that overshoot and undershoot the target profile,
as shown in
Fig.~\ref{fig:eps_dF}~(b). However,
due to the increase in the radial extent of the drug-loaded core that
occurs as $\epsilon$ increases, the penultimate overshoot in the efflux,
characterised by that which leads to its maximum value, has a greater
amplitude and occurs sooner in the drug-release process. Consequently,
the drug-release kinetics will be more prone to bursting when additional
drug molecules are loaded into the hydrogel and the drug-release period
is held constant. 

\begin{figure}
  \centering
  \includegraphics[width=0.48\textwidth]{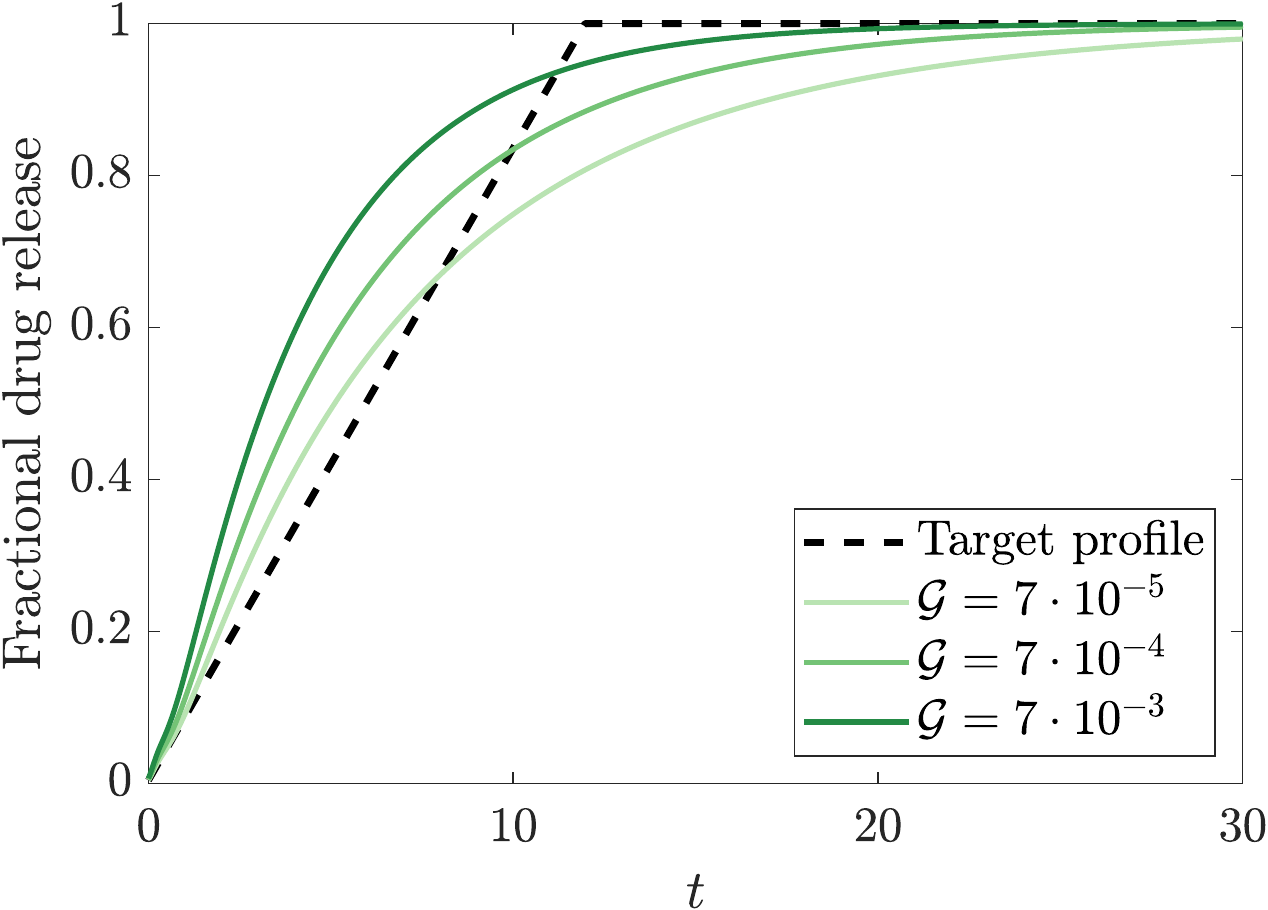}
  \caption{Optimised fractional drug-release profiles when the initial global
    drug fraction is fixed at $\epsilon = 0.1$. The other parameter values
    are $\chi = 0.5$ and $\D = 0.1$.}
  \label{fig:eps_G}
\end{figure}

When considering finite values of the global drug fraction
$\epsilon$, the hydrogel stiffness $\G$ can still be used to further tune the
drug-release profile. Specifically, stiffer and
softer gels accelerate and decelerate the release of drug
molecules, respectively; see Fig.~\ref{fig:eps_G}.
However, compared to the case of an infinitely
dilute drug, as shown in Fig.~\ref{fig:GeneralCaseFluxes}~(b)
for the same parameter set, we clearly
see that a finite value of $\epsilon$ generally leads to a more rapid
release of drug, particularly at small times. Thus, even the optimal
drug-release profiles exhibit stronger burst-like characteristics when
more drug is loaded into the hydrogel.
Tuning the
hydrogel stiffness may be an effective means of further mitigating the
burst effect in real-world applications, where the initial drug volume
is indeed finite, and overcome limitations with simply optimising over the
initial distribution of drug molecules.

\section{Conclusion}
\label{sec:conclusions}

The goal of this paper is to explore the possibility of optimising the
drug-release profile in hydrogel-based drug-delivery systems by
tuning the initial drug distribution. Thus, a model of
a spherical drug-loaded hydrogel that captures large elastic
deformations due to swelling is presented. By considering the limit of
a dilute drug, the equations for the hydrogel decouple from the equations
for drug transport, the latter of which become linear in the drug
concentration. Using this model, a theory for optimising the drug-release
profile via the
initial drug concentration is developed. The dilute-drug assumption
paves the way towards a fast numerical solver for the optimisation problem
by enabling the partial drug effluxes $f_i(t)$ used as basis functions
to be pre-computed.

For target drug effluxes that are piecewise-constant functions,
the optimal initial distribution of drug molecules generally
consists of a central drug-loaded packet or core with several isolated
packets near the
hydrogel boundary. The radial extent of the central core increases with the
amount of drug that is loaded into the hydrogel. Numerical simulations
reveal that the optimal initial drug concentrations are highly
effective at mitigating the burst effect and limiting the initial release
of drug into the surroundings. Moreover, the corresponding
drug-release profiles provide reasonable approximations to the target
profiles for all times.


The hydrogel stiffness provides an additional parameter that can be used
to tune the drug-release profile. However, the stiffness plays
a non-trivial role in the drug-release kinetics because it affects both
the drug diffusivity and the distance that drug molecules must travel to
reach the free surface. Softer hydrogels undergo a greater degree of
swelling; as a result, the drug diffusivity increases but so does the
distance to the free surface.
For the parameter range considered here, we find that the rate of drug
release 
is more strongly affected by variations in the radial extent of the hydrogel
than variations in the drug diffusivity. Consequently, softer hydrogels
lead to slower drug-release kinetics and are better suited for applications
that require the administration of drugs over long periods of time.
Moreover, the strong dependence of drug transport on the size of the
hydrogel highlights the importance of capturing finite deformations in
the model. 




The results presented here generate new questions about the impact of
viscoelasticity and network degradation on
the optimal loading of hydrogels. A reduction in
the elastic stress due to viscous rearrangement of the polymers
could trigger additional swelling and hence slow the
release of drug molecules.
Bulk degradation of the hydrogel, which could be captured by
a decreasing gel stiffness in time, could have similar consequences. 
By working within the dilute-drug limit, fast numerical methods can be
developed to optimise drug-release profiles calculated from  extended models.
A detailed study of the non-dilute limit would also be insightful by
capturing how the transport of drug molecules is affected by the mechanical
response of the hydrogel. Mathematical modelling will play a key role
in understanding these points and thus lead to finer control over the
delivery of drug payloads using hydrogel carriers.

\section*{Acknowledgements}

This research did not receive any specific grant from funding agencies in the public, commercial, or not-for-profit sectors.

\begin{appendix}
  \section{Existence and globality of local minima}
\label{app:optim}

In this appendix we prove that global minima exist and that local minima are
equivalent to global minima. To do so, it is helpful to recall that
a definition of convexity (for twice continuously differentiable functions) is that $ \forall \mathbf{v}$, $\mathbf{w} \in \R^{M}$ and $\alpha \in [0,1]$,
\begin{equation}
\label{eq:alphaineq}
     \alpha H(\mathbf{v}) + (1-\alpha)H(\textbf{w}) - H(\alpha \mathbf{v}  + (1-\alpha) \mathbf{w} )  \geq 0.
\end{equation}

We now consider the discrete optimisation problem given by
\eqref{eq:OverallOptim}--\eqref{dis:H}.
The feasible set for $\textbf{d}$ is bounded as $\textbf{d} \geq \textbf{0}$ and
\begin{equation}
    d_i = \int_0^{\infty} d_if_i(t)\, \d t \leq \int_0^{\infty} F(t)\, \d t =  \int_0^{\infty} A(t)\, \d t < \infty.
\end{equation}
Moreover, the feasible set is non-empty as all constraints are satisfied by,
for example,
\begin{equation}
  d_1 = \int_0^{\infty} A(t)\, \d t, \quad d_2, d_3, \ldots, d_M = 0.
\end{equation}
Since $H(\textbf{d})$ is continuous, the Boundedness Theorem shows that $H$ has a global minimum in this feasible set.

Using standard convex programming results~\cite{niculescu2006convex}, we can prove that constrained
local and global minima are equivalent. To do so, we consider a more
general convex programming problem given by
\begin{equation}
    \min \left\{H(\textbf{d}) : \textbf{L}\textbf{d} \leq \textbf{b} \right\},
\end{equation}
where $\textbf{L}$ is a matrix and $\textbf{b}$ a vector such that
$S := \{ \textbf{d} : \textbf{L}\textbf{d} \leq \textbf{b}\} \neq \emptyset$.
Note that the constraint $\textbf{L}\textbf{d} \leq \textbf{b}$ can
account for upper bounds on each $d_i$.
In this constrained optimisation problem, local minima can be defined by first considering the set of ``feasible directions'' $V(\textbf{d}^*)$ at a point $\textbf{d}^* \in S$. This is given by
\begin{equation}
    V(\textbf{d}^*) = \left\{ \boldsymbol{v} : \exists \beta > 0 \quad \text{such that} \quad \textbf{d}^* + \gamma \boldsymbol{v} \in S \quad \forall \gamma \in [0,\beta]\right\}.
\end{equation}
The vectors $\mathbf{v} \in V(\mathbf{d}^*)$ are the directions in which one could move a short distance from $\textbf{d}^*$ while remaining in the set
$S$. A local minimum $\textbf{d}^*$ therefore satisfies 
\begin{equation}
\label{eq:MinimumCond}
    \left.\td{}{\alpha}\left(H(\textbf{d}^* + \alpha \mathbf{v})\right)\right|_{\alpha = 0} \geq 0, \quad \forall \mathbf{v}\in V(\textbf{d}^*).
\end{equation}
To prove the necessity of \eqref{eq:MinimumCond},
suppose there exists a $\textbf{v} \in V(\textbf{d}^*)$ such that
\begin{equation}
\label{eq:derivcond}
    \left.\td{}{\alpha}\left(H(\textbf{d}^* + \alpha \mathbf{v})\right)\right|_{\alpha = 0} < 0. 
\end{equation}
Then, by continuity, (\ref{eq:derivcond}) implies that
\begin{equation}
\label{eq:derivcond2}
\exists \beta^* \text{ such that } \forall \delta \in [0,\beta^*],
\quad \left.\td{}{\alpha} \left(H(\textbf{d}^* + \alpha \mathbf{v})\right)\right|_{\alpha = \delta} < 0. 
\end{equation}
Then, for sufficiently small $\eta$,
\begin{equation}
    H(\textbf{d}^* + \eta \mathbf{v}) = H(\textbf{d}^*) + \int_0^{\eta} \left.\td{}{\alpha} \left(H(\textbf{d}^* + \alpha \mathbf{v})\right)\right|_{\alpha = \delta}\, \d \delta < H(\textbf{d}^*)
\end{equation}
and 
\begin{equation}
    \textbf{L}(\textbf{d}^* + \eta \mathbf{v}) \leq \textbf{b}
\end{equation}
so $\textbf{d}^*$ is not a global minimum.

Thus, it remains to show that all points that are not global minima also do not satisfy (\ref{eq:MinimumCond}). Suppose $\tilde{\textbf{d}}$ is a feasible point that is not a global minimum and that $\textbf{d}^*$ is a global minimum (which was proved to exist at the start of this section). Then,
\begin{equation}
\label{eq:wlogineq}
    H(\textbf{d}^*) < H(\tilde{\textbf{d}}).
\end{equation}
The convexity equation, (\ref{eq:alphaineq}), implies that
\begin{equation}
    \alpha H(\textbf{d}^*) + (1-\alpha)H(\tilde{\textbf{d}}) \geq H(\alpha \textbf{d}^*  + (1-\alpha) \tilde{\textbf{d}} ).
\end{equation}
Now, at $\alpha = 0$, the left- and right-hand sides are equal and so, for the inequality to hold,
\begin{equation}
    \td{}{\alpha} \left( \alpha H(\textbf{d}^*) + (1-\alpha)H(\tilde{\textbf{d}})\right) \geq \td{}{\alpha}\left(H(\alpha \textbf{d}^*  + (1-\alpha) \tilde{\textbf{d}} )\right) \space \space \space \quad \text{at}\quad \alpha = 0.
\end{equation}
Thus, by using (\ref{eq:wlogineq}), we obtain
\begin{equation}
\label{eq:LongIneq}
    0 >  H(\textbf{d}^*) - H(\tilde{\textbf{d}}) \geq \left.\td{}{\alpha} \left(H((\tilde{\textbf{d}} + \alpha(\textbf{d}^* - \tilde{\textbf{d}}))\right)\right|_{\alpha = 0}.
\end{equation}
Furthermore,
\begin{equation}
    \textbf{L}\left(\tilde{\textbf{d}} + \alpha(\textbf{d}^* - \tilde{\textbf{d}})\right)  = \textbf{L}\left(\alpha \textbf{d}^*  + (1-\alpha) \tilde{\textbf{d}}\right) \leq \textbf{b} \text{  } \quad \forall \alpha \in [0,1],
\end{equation}
by feasibility of $\tilde{\textbf{d}}$ and $\textbf{d}^*$ which means that $\textbf{v} := \textbf{d}^* - \textbf{d}$ is a feasible direction (by defining $\beta := 1$). Together with (\ref{eq:LongIneq}) this shows that $\tilde{\textbf{d}}$ is not a local minimum as required and so local and global minima are equivalent for this problem.

\end{appendix}

\bibliography{refs}
\bibliographystyle{unsrtnat}

\end{document}